\documentclass[12pt]{amsart}
\usepackage{amsmath}
\usepackage{amsthm}
\usepackage[T1]{fontenc}
\usepackage[sc]{mathpazo}
\usepackage{enumerate}

\usepackage{etoolbox}
\patchcmd{\section}{\scshape}{\bfseries}{}{}

\newenvironment{nouppercase}{%
  \renewcommand{\uppercasenonmath}[1]{}}{}

\newcommand{\C}{{\mathbb{C}}}
\newcommand{\R}{{\mathbb{R}}}
\newcommand{\N}{{\mathbb{N}}}
\newcommand{\Z}{{\mathbb{Z}}}

\newcommand{\conj}[1]{\overline{#1}}

\newcommand{\pp}{\widetilde p}
\newcommand{\PP}{\widetilde P}

\DeclareMathOperator{\Log}{Log}

\newtheorem{theorem}{Theorem}[section]
\newtheorem{proposition}[theorem]{Proposition}
\newtheorem{lemma}[theorem]{Lemma}

\newtheorem{corollary}[theorem]{Corollary}
\theoremstyle{definition}

\theoremstyle{remark}
\newtheorem{remark}[theorem]{Remark}
\newtheorem*{remark*}{Remark}
\theoremstyle{example}

\begin{document}
\title[~]
{Characterization of polynomials whose large powers \\
have fully positive coefficients}
\author[~]
{Colin Tan and Wing-Keung To}

\address{Colin Tan,
Engineering Systems and Design, Singapore University of Technology and Design,
8 Somapah Road, Singapore 487372
}
\email{colinwytan@protonmail.com}

\address{Wing-Keung To,
Department of Mathematics, National University of Singapore, Block
S17, 10 Lower Kent Ridge Road, Singapore 119076}
\email{mattowk@nus.edu.sg}

\thanks{Wing-Keung To was partially supported by the Singapore Ministry of Education Academic Research  Fund Tier 1 grant
R-146-000-254-114.
}

\keywords{
polynomials, positive coefficients, projective toric manifolds 
}
\subjclass[2010]{
26C05, 14M25, 32L15
}
\begin{nouppercase}
\maketitle
\end{nouppercase}
\numberwithin{equation}{section}

\begin{abstract}
We give a criterion which characterizes a 
real multi-variate Laurent polynomial with full-dimensional smooth Newton polytope to have the property that
all sufficiently large powers of the polynomial have fully positive coefficients.  Here a Laurent polynomial is said to have fully positive coefficients if the coefficients of its monomial terms indexed by the lattice points of its Newton polytope are all positive.  Our result generalizes 
an earlier result of the authors, which corresponds to the special case when the Newton polytope of the Laurent polynomial is a translate of a standard simplex.  The result also generalizes a result of De Angelis, which corresponds to the special case of univariate polynomials.
As an application, we also give a characterization of certain 
polynomial spectral radius functions of the defining matrix functions of Markov chains.
\end{abstract}

\section{Introduction and main results} \label{sec: intro}

Positivity conditions for polynomials with real coefficients play a key role in several branches of  mathematics, such as real algebraic geometry, convex geometry, probability theory and 
optimization, and have been widely studied (see e.g. \cite{CD97, Deangelis94, 
Deangelis03, Han86, MT93, Polya28, PR09, Reznick95, ST17, TT18, TY06}  and the references therein).  An interesting and important class of polynomials are those whose coefficients are positive. 
				
\medskip
A polynomial $p$ of degree $d$ in $n$ variables, upon homogenizing, gives rise to a homogeneous polynomial 
$\pp$ of the same degree $d$ in $n+1$ variables.  In a recent work \cite{TT18}, the authors obtained a 
characterization of those polynomials $p$ in $n$ variables possessing the property that 
	$p^m$ has all positive coefficients (as an inhomogeneous polynomial of degree $md$) for all sufficiently large $m$ in terms of certain positivity conditions on the associated homogeneous 
								polynomial $\pp$.  The work generalizes an earlier result of De Angelis \cite[Theorem 6.6]{Deangelis03}, which corresponds to the case when $n=1$.  We remark that the Newton polytope  of any $p$ having the above property is necessarily the standard simplex $\Delta_{n,d}=:\{(x_1,\cdots,x_n)\,\big|\, x_i\geq 0, ~i=1,\cdots,n, \text{ and }x_1+\cdots +x_n\leq d\}$ in $\R^n$ of length $d$.  As such, it is interesting and natural to ask whether similar result holds in the more general setting when the Newton polytope of $p$ need not be a standard simplex.
																				
\medskip							
In this paper, we generalize the afore-mentioned result 
to the case when the Newton polytope of the polynomial is a full dimensional smooth lattice polytope.	Let $n \ge 1$, and let $\displaystyle p= \sum_{ m\in\Z^n } c_{m} x^m\in \R[x_1^{\pm 1},\dots,x_n^{\pm 1}]$ be a
Laurent polynomial with real coefficients in the $n$ variables $x_1,\cdots,x_n$.  Then $p$ is said to 
		{\emph{have fully positive coefficients}} if $c_m> 0$ for all $m\subset \Phi \cap\Z^n$,
		where $\Phi$ denotes the Newton polytope of $p$.  
Here $m= (m_1,\ldots, m_n)\in\Z^n$ is a multi-index
			and
				$x^{m} = x_1^{m_1} x_2^{m_2} \cdots x_n^{m_n}$. 
Let $\Sigma$ be the normal fan of $\Phi$.  When $\Phi$ is an $n$-dimensional 
smooth lattice polytope, the toric variety $X_\Sigma$ associated to $\Sigma$ is an $n$-dimensional projective manifold. For each $0\leq i\leq  n$ (and each cone $\sigma\in\Sigma$), we denote the collection of $i$-dimensional cones (resp. faces) of $\Sigma$ (resp. $\sigma$) by $\Sigma(i)$ (resp. $\sigma(i)$).  Let 
$\R_+:=\{x\in\R\,\big|\,x\geq 0\}$, and consider the
circle group $U(1) := \{e^{i\theta} \, \big| \, \theta\in \R\}\subset \C^*:=\C\setminus\{0\}$.
Denote by $\C^{\Sigma(1)}$ (resp. $(\C^*)^{\Sigma(1)}$, $\R^{\Sigma(1)}$, $\R_+^{\Sigma(1)}$, $
U(1)^{\Sigma(1)}$) the Cartesian product of copies of 
$\C$ (resp. $\C^*$, $\R$, $\R_+$, 
$U(1)$) indexed by $\rho\in \Sigma(1)$, etc.  In particular, $\R_+^{\Sigma(1)}$ is the closed positive orthant of $\R^{\Sigma(1)}$.
Let $\displaystyle Z(\Sigma):=
\{z=(z_\rho)_{\rho\in\Sigma(1)}\in
\C^{\Sigma(1)}\,\big|\,\prod_{\rho\in \Sigma(1)\setminus\sigma(1)}z_\rho=0\text{ for all }\sigma\in\Sigma(n)\}$.  It is known that $X_\Sigma=(\C^{\Sigma(1)}\setminus Z(\Sigma) )/G$
for some subgroup
$G\subset (\C^*)^{\Sigma(1)}$ which is naturally isomorphic to $\text{Hom}_{\Z}(\text{Cl}(X_\Sigma),\C^*)$ (see \cite{CLS11} or \cite{Cox95}).  Here $(\C^*)^{\Sigma(1)}$ acts on $\C^{\Sigma(1)}$ by coordinatewise multiplication, and
$\text{Cl}(X_\Sigma)$ denotes the Weil divisor class group of $X_\Sigma$.
Consider the polynomial ring
$\C[z_\rho\,\big|\,\rho\in \Sigma(1)]$, and let $\pp\in \C[z_\rho\,\big|\,\rho\in \Sigma(1)] $ denote the 
$\Phi$-homogenization of $p$ (see e.g. \cite[\S 5.4]{CLS11}).  
For each $\sigma\in\Sigma(n)$, we let $e^{(\sigma)}=(e^{(\sigma)}_\rho)_{\rho\in\Sigma(1)}
\in \C^{\Sigma(1)}$ be the point given by 
\begin{equation}\label{eq: defesigma}
e^{(\sigma)}_\rho=\begin{cases}1   \quad \text{if }\rho\in\Sigma(1)\setminus\sigma(1),
\\
0\quad \text{if }\rho\in\sigma(1).
\end{cases}
\end{equation}
Let
$(G\cap U(1)^{\Sigma(1)})\cdot \R_+^{\Sigma(1)}:=\{ g\cdot x\,\big|\, g\in G\cap U(1)^{\Sigma(1)}\text{ and }x\in\R_+^{\Sigma(1)}\}$, where $g\cdot x$ denotes the coordinatewise product of the tuples $g$ and $x$.
For each $\rho_o\in \Sigma(1)$, we denote the associated facet of $\R_+^{\Sigma(1)}$  given by	
$
F_{\rho_o}(\R_+^{\Sigma(1)}) := \{ x=(x_\rho)_{\rho\in \Sigma(1)}\in \R_+^{\Sigma(1)} \, | \,  x_{\rho_o}= 0\}$. 

\medskip
Notation as above.  
Our main result in this paper is the following: 
				
\begin{theorem} \label{thm: mainThm}
Let $p \in \R[x_1^{\pm 1},\dots,x_n^{\pm 1}]$ be a Laurent polynomial whose Newton polytope $\Phi$ is an
$n$-dimensional smooth lattice polytope. Let $\pp$ be the $\Phi$-homogenization of $p$.  
The following two statements are equivalent:
\begin{enumerate}[(a)]
\item \label{thm: Condition}
$\pp$ satisfies the following three conditions:
\begin{description}
\item[(Pos1)] \label{Pos1} $\pp(e^{(\sigma)})>0$ for all $\sigma\in \Sigma(n)$.
\item[(Pos2)] \label{Pos2} For all $\rho\in \Sigma(1)$,
								 $\displaystyle\frac{\partial \pp}{ \partial x_\rho} (x) > 0$ 
  											for all  $ x \in F_{\rho}(\R_+^{\Sigma(1)}) \setminus
  											 (Z(\Sigma)\cap F_{\rho}(\R_+^{\Sigma(1)}) )$.
\item[(Pos3)]  \label{Pos3} $|\pp(z)| < \pp((|z_\rho|)_{\rho\in \Sigma(1)})$ for all 
$z=(z_\rho)_{\rho\in\Sigma(1)} \in \C^{\Sigma(1)} \setminus (Z(\Sigma)\cup (G\cap U(1)^{\Sigma(1)})\cdot \R_+^{\Sigma(1)})$.
\end{description}

\item \label{thm: EventualPos} There exists $k_o > 0$ such that for each integer $k \ge  k_o$, 
								\, $p^k$ has fully positive coefficients.
\end{enumerate}
\end{theorem}								
								
We refer the reader to Section \ref{sec: positivity} for the definitions and/or detailed discussion of the various terms in Theorem \ref{thm: mainThm}.

\medskip
Compared to \cite{TT18}, a main new ingredient in this paper is that we exploit extensively the toric geometry that underlies Newton polytope $\Phi$ of $p$.  
 In particular, $p$ gives rise to a Hermitian algebraic function on a very ample holomorphic line bundle over 
the projective toric manifold $X_\Sigma$.  
The bulk of our proof of the implication $\eqref{thm: Condition} \implies \eqref{thm: EventualPos}$ 
consists of showing that the Hermitian algebraic function satisfies the sufficient conditions of a Hermitian Positivstellensatz of Catlin-D'Angelo 
\cite{CD99} (see also Theorem \ref{TheoremCD99} below), enabling us to apply the latter result.  

\medskip
An interesting example of polynomials satisfying the three positivity conditions in (a) and modified from 
D'Angelo-Varolin \cite[Theorem 3]{DV04}
is given by
\begin{align}\label{eq: plambda}
p_{\lambda_1,\lambda_2}(x_1, x_2) :&= \prod_{i=1,2} \big[(1+x_i)^{2\ell} - \lambda_i x_i^\ell\big]\\
\nonumber
\quad\text{with  }
\ell \ge 2&\text{ and }
\binom{2\ell}{\ell}<\lambda_i < 2^{2\ell - 1},~i=1,2,
\end{align}
where the Newton polytope $\Phi$ of $p_{\lambda_1,\lambda_2}$ is the square $[0,2\ell]^2$ (and not a simplex in $\R^2$), and 
 the coefficients of $x_1^\ell x_2^j$ and $x_1^j x_2^\ell$, $j\in \{0,1,\cdots,\ell-1,\ell+1\cdots,2\ell\}$,
are negative.  Nonetheless the $\Phi$-homogenization of $p_{\lambda_1,\lambda_2}$ satisfies the three 
positivity conditions in $\eqref{thm: Condition}$ (we will skip the verification which is similar to the calculations in 
\cite{DV04}), and thus Theorem \ref{thm: mainThm} is applicable to $p_{\lambda_1,\lambda_2}$.

\medskip
The three positivity conditions in $\eqref{thm: Condition}$ are independent of each other.  We refer the reader to \cite[Section 1]{TT18} for examples of polynomials (in the simplex case) which satisfy two of the conditions but do not satisfy all three of them.

\medskip
In view of Theorem \ref{thm: mainThm}, it is natural to ask for a similar characterization of 
polynomials whose large powers have \lq fully nonnegative coefficients\rq.   Another natural question that arises is whether Theorem \ref{thm: mainThm} generalizes to the case when the Newton polytope $\Phi$ is not smooth.  The method in this paper does not appear to generalize readily to handle such cases, and new ideas will be 
needed.
To glimpse the intricacy of the first question, we mention that a
limiting case of the family of polynomials in \eqref{eq: plambda}, namely $p_{2^{2\ell - 1},2^{2\ell - 1}}$ (so that
$\lambda_1=\lambda_2= 2^{2\ell - 1}$), satisfies (Pos1), (Pos2) and 
a weaker version of (Pos3) (with \lq$<$\rq~there replaced by  \lq$\leq$\rq), but one easily checks that 
all of its powers have some negative coefficients.  As for the non-smooth case, one likely needs some kind of Hermitian Positivstellensatz of Catlin-D'Angelo 
\cite{CD99} for singular varieties, which to the authors' knowledge, has not been established yet.



\medskip
An interesting question in the study of Markov chains is to characterize those polynomials $q$ 
        for which there exists an irreducible (resp. aperiodic) Markov chain whose defining matrix (or equivalently, an irreducible (resp. aperiodic) square matrix over $\Z_+[x_1,\ldots, x_\ell]$) has $q$ as its 
        spectral radius function (see e.g. \cite{Deangelis94, Deangelis942}).  Here, $\Z_+:=\{k\in\Z\,\big|\, k\geq 0\}$,  and $\Z_+[x_1,\ldots, x_\ell]$ denotes the 
semiring of polynomials in $x_1,\ldots, x_\ell$ with coefficients in $\Z_+$.  This spectral radius function
         is an important invariant in the study of Markov shifts
                (see e.g. \cite{MT93}).  
Denote the spectral radius function of an {\it irreducible } (resp. {\it aperiodic}) 
square matrix $A$
over $\Z_+[x_1,\ldots, x_\ell]$
by $\beta_A=\beta_A(x_1,\ldots, x_\ell)$.
As an application of Theorem \ref{thm: mainThm}, we have

\begin{corollary} \label{cor: betaApplication}
Let $p \in \Z[x_1^{\pm 1},\dots,x_n^{\pm 1}]$ be a Laurent polynomial whose Newton polytope $\Phi$ is an $n$-dimensional smooth lattice polytope, and such that the $\Phi$-homogenization $\pp\in\Z[ z_\rho\,\big|\,\rho\in\Sigma(1)]$ of $p$ satisfies (Pos1) and (Pos2).
The following statements are equivalent:
\begin{enumerate}[(i)]
\item \label{cor: pos3}$\pp$ satisfies (Pos3).
\item \label{cor: irr} $\pp = \beta_A$ for some irreducible square matrix $A$ over $\Z_+[ z_\rho\,\big|\,\rho\in\Sigma(1)]$.
\item \label{cor: aper} $\pp = \beta_A$ for some aperiodic square matrix $A$ over 
$\Z_+[ z_\rho\,\big|\,\rho\in\Sigma(1)]$.
\end{enumerate}
\end{corollary}

Similar to Theorem \ref{thm: mainThm}, Corollary \ref{cor: betaApplication} generalizes an earlier result of the authors 
\cite[Corollary 1.2]{TT18} which corresponds to the case when $\Phi$ is a standard simplex, as well as a result of De Angelis \cite[Theorem 6.7]{Deangelis942} which corresponds to the case when $n=1$.  
We refer the reader to 
\cite[Section 5]{TT18} for a convenient recollection of the definition of an \lq{\it irreducible } (resp. {\it aperiodic}) 
square matrix 
over $\Z_+[ z_\rho\,\big|\,\rho\in\Sigma(1)]$\rq, and we remark that 
the interpretation of Corollary \ref{cor: betaApplication} in terms of Markov chains for the standard simplex case as given in \cite[Section 5]{TT18} also holds in the present more general setting.

\medskip
De Angelis' Positivstellensatz \cite[Theorem 6.6]{Deangelis03} has been applied by Bergweiler-Eremenko \cite{BE15}
	to study the distribution of zeros of polynomials with positive coefficients 
		(see also \cite{EF}).		
As a generalization of the Positivstellensatz of De Angelis,
Theorem \ref{thm: mainThm}
  may also have similar applications,
  	which we will not pursue here.

\medskip
The organization of this paper is as follows. In Section \ref{sec: CDThm}, we recall 
some background material
on
Hermitian algebraic functions on holomorphic line bundles. 
In Section \ref{sec: positivity}, we recall the
toric geometry associated to a real Laurent polynomial, and 
relate some positivity properties of the
polynomial with those of its $\Phi$-homogenization.
In Section \ref{sec: proofOfTheorem}, we give the proof of Theorem \ref{thm: mainThm}.  
In Section \ref{sec: Application}, we give the deduction of Corollary \ref{cor: betaApplication}.
 
 \medskip
\textbf{Acknowledgements.}
The authors would like to acknowledge 
	John P. D'Angelo,
	Valerio De Angelis and
	David Handelman for sharing their work and for helpful discussions.

\section{Hermitian algebraic functions and Catlin-D'Angelo's Positivstellensatz}
	\label{sec: CDThm}

In this section, we recall some background material regarding Hermitian algebraic functions
which is mostly taken from \cite{CD97, CD99, DAngelo02, DV04, Var08}.  
 As such, we will skip their proofs here and refer the reader to these references for their proofs.

\medskip
Let
$X$ be an $n$-dimensional compact complex manifold, and let $F$ be a holomorphic line bundle over $X$ with its projection map denoted by $\pi:F\to X$.    The dual holomorphic line bundle of $F$ is denoted by $F^*$, and the complex conjugate manifold (resp. bundle) of $M$ (resp. $F$) is denoted by
$\overline{X}$ (resp. $\overline{F}$), etc.
 Let $\pi_1:X\times \overline{X}\to X$ and $\pi_2:X\times\overline{X}\to \overline{X}$ denote the projection maps onto the first and second factor respectively, and consider the holomorphic line bundle $\pi_1^*F\otimes \pi_2^*\overline{F}$ over the complex manifold $X\times \overline{X}$, whose fiber at a point $(x,\overline{y})\in X\times \overline{X}$ is naturally isomorphic to $F_x\otimes \overline{F_y}$.  Here $F_x:=\pi^{-1}(x)$ denotes the fiber of $F$ at the point $x\in X$, etc. 
 Following \cite{Var08}, a {\it Hermitian algebraic function $Q$ on $F$} is defined as a global holomorphic section of $\pi_1^*F\otimes \pi_2^*\overline{F}$ (i.e., $Q\in H^0(X\otimes  \overline{X}, \pi_1^*F\otimes \pi_2^*\overline{F})$) satisfying the condition
\begin{equation}\label{eq:  HermitianQ}
Q(x,\overline{y})=\overline{Q(y,\overline{x})}\in F_x\otimes\overline{F_y}\quad\text{for all }x,y
\in X
\end{equation}
 (see also \cite{Var08} for an alternative definition of $Q$ as a function on $F^*\times \overline{F^*}$
satisfying analogous conditions).
One easily sees that with respect to any basis $\{s^\alpha\}$ of $H^0(X, F)$,
there exists a corresponding Hermitian matrix $\big( C_{\alpha\overline{\beta}}\big)$ such that,
for all $x,y\in X$, one has
$\displaystyle
Q(x,\overline{y})=\sum_{\alpha,\beta}C_{\alpha\overline{\beta}}s^\alpha(x)\overline{s^\beta(y)}
$.  
One says
 that $Q$ is {\it a maximal sum of Hermitian squares} if the Hermitian matrix
$\big(C_{\alpha\overline{\beta}}\big)$ with respect to one (and hence any) basis of $H^0(X, F)$ is positive definite, or equivalently, there exists a basis $\{t^\alpha\}$ of $H^0(X, F)$ such that $Q(x,\overline{x})=\sum_\alpha t^\alpha(x)\overline{t^\alpha(x)}$ for all $x\in X$.
Also, the Hermitian algebraic function $Q$ is said to be {\it positive} if 
\begin{equation}\label{eq: positiveQ}
Q(v,\overline{v})>0\quad \text{for all }0\neq v\in F^*.
\end{equation}  Here and henceforth, with slight abuse of notation,  $Q(v,\overline{v})$ denotes the obvious evaluation induced by the natural pairing between $F$ and $F^*$.  If $Q$ is positive, then $Q$ induces a Hermitian metric $h_Q$ on $F^*$ given by $h_Q(v,w)=Q(v,\overline{w})$ for $v,w\in F^*_x$, $x\in X$.  (In \cite{CD99}, such a Hermitian metric arising from a positive Hermitian algebraic function is called a {\it globalizable metric}.)  We recall that the curvature form $\Theta_{h_Q}$ of the Hermitian metric $h_Q$ is the $(1,1)$-form on $X$ given locally as follows:   On
any open subset $U$ of $X$ and for any local non-vanishing holomorphic section $s$ of $F^*\big|_U$, one has $\Theta_{h_Q}\big|_U=-\sqrt{-1}\partial\overline\partial\log h_Q(s,\overline{s})=-\sqrt{-1}\partial\overline\partial\log Q(s,\overline{s})$.  Following \cite{Var08} again (and with origin in \cite{CD99}), a positive Hermitian algebraic function $Q$ on $X$ is said to satisfy the {\it strong global Cauchy-Schwarz (SGCS) inequality} if
\begin{equation}\label{eq: SGCSineq}
|Q(v,\overline{w})|^2<Q(v,\overline v)\cdot Q(w,\overline w) 
\end{equation}
\text{for all non-zero }$v,w\in F^*$\text{ such that }$\pi(v)\neq \pi(w)$.
Note that one always has $|Q(v,\overline{w})|^2=Q(v,\overline v)\cdot Q(w,\overline w)$ whenever $\pi(v)= \pi(w)$.

\medskip
We recall the following result of Catlin-D'Angelo:

\begin{theorem}  [\cite{CD99}]\label{TheoremCD99}
Let $L $ and $E$ be  holomorphic line bundles over an $n$-dimensional compact complex manifold $X$.
Suppose $R$ and $Q$ are positive Hermitian algebraic functions on $L$ and $E$ respectively, such that 
\begin{enumerate}[(i)]
\item \label{DA: SGCS}$R$ satisfies the SGCS inequality, and 
\item \label{DA: curv}the curvature $(1,1)$-form $\Theta_{h_R}$ on $X$ is negative definite.   
\end{enumerate}
Then there exists $m_o\in \mathbb N$ such that for each integer
$m\geq m_o$, the tensor product $R^mQ$ is a maximal sum of Hermitian squares (on the line bundle $L^{\otimes m}\otimes E$).
	\end{theorem} 

We remark that  Catlin-D'Angelo obtained the above theorem by proving the positive-definiteness of certain associated integral operators on $X$.

\medskip
\section{Laurent polynomials and associated toric geometry} 
\label{sec: positivity}

In this section, we recall some background material on the toric geometry that underlies Laurent polynomials, which is mostly taken from \cite{CLS11, Cox95, Ful93}.  We will also relate the positivity conditions on Laurent polynomials with those on associated Hermitian algebraic functions on holomorphic line bundles over the underlying 
projective toric manifold as well as those of associated polynomials on the total coordinate ring.

\medskip
Throughout this section, we fix a positive integer $n\geq 1$.  Let $\displaystyle p= \sum_{ m\in\Z^n } c_{m} x^m\in \R[x_1^{\pm 1},\dots,x_n^{\pm 1}]$ be a
Laurent polynomial with real coefficients in the $n$ variables $x_1,\cdots,x_n$.  Here $m= (m_1,\ldots, m_n)\in\Z^n$ is a multi-index
			and
				$x^{m} = x_1^{m_1} x_2^{m_2} \cdots x_n^{m_n}$.  Consider the finite set
$\Log(p):=\{m\in\Z^n\,\big|\,\, c_{m} \neq 0 \}$.  Then the {\emph{Newton polytope}} $\Phi=\Phi(p)$ of $p$ is defined as the convex hull
	of 
	$\text{Log}(p)$ in $\R^n$.  Following \cite{CLS11, Ful93}, we write 
	\begin{equation}\label{eq: MZn}
	 M:=\Z^n,\quad\text{so that}\quad\Log(p)\subset M,
	\end{equation} 
	and $\Phi$ is a lattice polytope in $M_\R:=M\otimes_{\Z}\R\cong \R^n$.  The dual lattice of $M$ is denoted by $N:=\text{Hom}_\Z(M,\Z)\cong \Z^n$, and one writes $N_\R:=N\otimes_{\Z}\R\cong \R^n$. 
A {supporting affine hyperplane} of $\Phi$ is an affine hyperplane $H$ in $M_R$ such that $\Phi$ lies entirely on one side of $H$, and a {\it face} of $\Phi$ is a non-empty set of the form $H\cap \Phi$ for some supporting affine hyperplane $H$.  Throughout this section, we will assume that $\Phi$ is of affine dimension $n$.  As such, the {\it vertices}, {\it edges} and {\it facets} of $\Phi$ are its faces of affine dimension $0$, $1$ and $n-1$ respectively.  For each $0\leq i\leq n$, we denote the set of $i$-dimensional faces of 
$\Phi$ by $\Phi(i)$. 
For each vertex $v$ of $\Phi$, we denote the set of edges of $\Phi$ containing $v$ by $\Phi(1)_v$; and for each $E\in \Phi(1)_v$, we 
denote 
by $w_{E,v}$ the point in $(E\cap M)\setminus \{v\}$ which is adjacent to $v$.  Throughout this section, we will assume that the lattice polytope $\Phi$ is {\it smooth}, which means that for each vertex $v$ of $\Phi$, the set $\{w_{E,v}-v\,\big|\,E\in \Phi(1)_v \}$ forms a $\Z$-basis of $M$.   It is easy to see that for each facet $F\in \Phi(n-1)$, there exists a unique minimal inward pointing normal $u_F\in N$ and a unique number $a_F\in \Z$ such that the unique supporting affine hyperplane $H_F$ of $F$ is given by
\begin{equation}\label{eq: supporthyperplane}
H_{F}:=\{m\in M_\R\,\big|\, \langle m,u_F\rangle=-a_F\},
\end{equation}
where $ \langle ~,~\rangle:M_{\R}\times N_{\R}\to \R$ is the natural pairing between $M_\R$ and $N_\R$.
Furthermore, one has the unique facet presentation of $\Phi$ given by
\begin{equation}\label{eq: facetpresentation}
\Phi=\{m\in M_\R\,\big|\, \langle m,u_F\rangle\geq -a_F\text{ for all }F\in \Phi(n-1)\}
\end{equation}
 (see e.g. \cite[\S 2.2]{CLS11}).  For a finite subset $S\subset N_\R$, we denote the cone spanned by $S$ in $N_\R$ by $\text{Cone}(S):=\{\sum_{w\in S}\lambda_w w\,\big|\, \lambda_w\geq 0\}$.  
 For each face $Q$ of $\Phi$, we denote the associated cone in $N_\R$ given by 
 $\sigma_Q:=\text{Cone}(\{u_F\,\big|\,F\in \Phi(n-1),~F\supset Q\})$, where each $u_F$ is as in \eqref{eq: supporthyperplane}.
 Then the {\it normal fan} $\Sigma=\Sigma(\Phi)$ of $\Phi$ is given by 
 \begin{equation}\label{eq: normalfan}
\Sigma:=\{\sigma_Q\,\big|\, Q \text{ is a face of }\Phi\}.
 \end{equation}
Note that $\Sigma$ is a fan in the sense that any face of a cone in $ \Sigma$ is itself a cone in $ \Sigma$, and 
the intersection of two cones in $\Sigma$ is a face of both cones. 
For each $0\leq i\leq n$, we denote the set of $i$-dimensional cones in $\Sigma$ by $\Sigma(i)$ and the $i$-dimensional faces of a cone $\sigma\in\Sigma$ by $\sigma(i)$.  Then one easily sees that there is a bijection between $\Phi(i)$ and $\Sigma(n-i)$ given by
 \begin{equation}\label{eq: coneface}
 Q\in \Phi(i)~~\longleftrightarrow ~~\sigma_Q\in \Sigma(n-i).
 \end{equation}
For each $\rho\in \Sigma(1)$, we write $u_\rho:=u_F$ and $a_\rho:=a_F$, where $F\in\Phi(n-1)$ is the facet that corresponds to $\rho$ under \eqref{eq: coneface}.  The supporting hyperplane $H_F$ of the above $F$ will also be denoted by $H_\rho$.  Then we may rewrite 
\eqref{eq: facetpresentation} as
\begin{equation}\label{eq: dualfacetpresentation}
\Phi=\{m\in M_\R\,\big|\, \langle m,u_\rho\rangle\geq -a_\rho\text{ for all }\rho\in \Sigma(1)\}.
\end{equation}
Associated to the fan $\Sigma$ is a toric variety $X_\Sigma$ containing the torus $T_N\cong  \text{Hom}_\Z(M,\C^*) \cong N\otimes_\Z \C^*\cong(\C^*)^n$ as a Zariski open subset and such that the natural action of $T_N$ on itself extends to a morphism $T_N\times X_\Sigma\to X_\Sigma$.  Since $\Phi$ is an $n$-dimensional smooth lattice polytope, it follows that the support of $\Sigma$ is $N_\R$, $\Sigma$ is a smooth fan (i.e., each cone $\sigma$ in $\Sigma$ is a smooth cone in the sense that the minimal generators of the lattice points of the rays in $\sigma(1)$ form 
a part of a 
$\Z$-basis of $N$), 
and as a consequence, $X_\Sigma$ is an $n$-dimensional (smooth) projective manifold (cf. e.g. \cite[\S 2.4, \S 3.1]{CLS11}).  
From the Orbit-Cone Correspondence, for each $i$-dimensional cone $\sigma\in\Sigma$, one has an associated $(n-i)$-dimensional $T_N$-orbit in $X_\Sigma$ which will be denoted by $O(\sigma)$
 (see e.g. \cite[\S 3.2]{CLS11}).  Furthermore, for each $\rho\in \Sigma(1)$, the closure of $O(\rho)$ in 
 $X_\Sigma$ is an $T_N$-invariant Cartier divisor, which will be denoted by $D_\rho$.  We denote the 
 holomorphic line bundle over $X_\Sigma$ associated to the divisor 
 $\sum_{\rho\in\Sigma(1)}a_\rho D_\rho$ by
 \begin{equation}\label{eq:  defL}
L:=\big[\sum_{\rho\in\Sigma(1)}a_\rho D_\rho\big],
 \end{equation}
where each $a_\rho$ is as in \eqref{eq: dualfacetpresentation}.   Since $\Phi$ is a smooth polytope, one knows that $L$ is very ample (see e.g. \cite[Proposition 6.1.10 and Theorem 6.1.15]{CLS11}).  Each $m=(m_1,\cdots,m_n)\in M$ gives rise to a character $\chi^m$ on $T_N\cong (\C^*)^n$ (with the isomorphism induced from \eqref{eq: MZn}) given by $\chi^m(t_1,\cdots,t_n)=t_1^{m_1}\cdots t_n^{m_n}$ for each $t=(t_1\cdots,t_n)\in T_N$.  It is known that for each $m\in M$ and each $\rho\in \Sigma(1)$, the vanishing order of $\chi^m$ along $D_\rho$ is given by $\langle m,u_\rho\rangle$.  Together with \eqref{eq: dualfacetpresentation}, it follows that for each 
 $m\in \Phi\cap M$, $\chi^m$ may be regarded as a global holomorphic section of $L$.  In fact, it is well-known that under such identification, one has the following isomorphism given by 
 \begin{equation}\label{eq: globalsectionisomorphism}
 H^0(X_\Sigma, L)\cong \bigoplus_{m\in\Phi\cap M}\C\cdot\chi^m
 \end{equation}
 (see e.g. \cite[Proposition 4.3.3]{CLS11}).  Next we associate to $p$ (via 
 \eqref{eq: globalsectionisomorphism}) the following $\widehat p\in  H^0(X_\Sigma, L)$ and Hermitian algebraic function $\widehat P\in  H^0(X_\Sigma\times \overline{X_\Sigma}, \pi_1^*L\otimes \pi_2^*\overline{L})$ given by
  \begin{align}\label{eq: Pdef1}
\widehat p(x)&:= \sum_{m\in \Phi\cap M }c_m\cdot\chi^m(x)\quad\text{and}\\
\label{eq: Pdef2}
\widehat P(x,\overline y)&:= \sum_{m\in\Phi\cap M }c_m\cdot\chi^m(x)\cdot \overline{\chi^m(y)}\quad\text{for }x,y\in X_\Sigma,
 \end{align}
where $\pi_i$ denotes the projection of $X_\Sigma\times \overline{X_\Sigma}$ onto the $i$-th factor, $i=1,2$.  Here and henceforth, $c_m$ is understood to be $0$ if $m\in(\Phi\cap M)\setminus \text{Log}(p) $ (note that $\text{Log}(p)\subset \Phi\cap M$).  
We remark that $\widehat P$ satisfies 
\eqref{eq:  HermitianQ} since $c_m\in \R$ for all $m\in\Phi\cap M $.  Recall also that $p$ is said to {\it have fully positive coefficients} if 
$c_m>0$ for all $m\in \Phi\cap M $.  First we have

\begin{proposition}\label{prop: AssocBihomIsSOS}
$p$ has fully positive coefficients
    if and only if
      $\widehat P$
            is a maximal sum of Hermitian squares.
\end{proposition}

\begin{proof}  As indicated in \eqref{eq: globalsectionisomorphism}, a basis of $H^0(X_\Sigma, L)$ is given by
$\{\chi^m\,\big|\, m\in\Phi\cap M\}$.
With respect to this basis, it follows from 
    \eqref{eq: Pdef2} that the square matrix associated to $\widehat P$ is given by 
    the real diagonal matrix $C:=\mathop{\mathrm{diag}}(c_m)_{m\in\Phi\cap M}$.  Then, as remarked in 
    Section \ref{sec: CDThm},
$\widehat P$
            is a maximal sum of Hermitian squares 
    if and only if $C$ is a positive definite matrix.  
    In turn, the latter condition holds if and only if $c_m>0$ for all $m\in\Phi\cap M$, or equivalently,
    $p$ has fully positive coefficients.
\end{proof} 

\medskip
As defined in Section \ref{sec: intro}, we have 
$\R_+=\{x\in\R\,\big|\,x\geq 0\}$, 
$U(1) = \{e^{i\theta} \, \big| \, \theta\in \R\}\subset \C^*$,
and
$\C^{\Sigma(1)}$ (resp. $(\C^*)^{\Sigma(1)}$, $\R^{\Sigma(1)}$, $\R_+^{\Sigma(1)}$, $
U(1)^{\Sigma(1)}$) is the Cartesian product of copies of 
$\C$ (resp. $\C^*$, $\R$, $\R_+$, 
$U(1)$) indexed by $\rho\in \Sigma(1)$, etc.  Following \cite[\S 5.1]{CLS11}, the {\it total coordinate ring}
of $X_\Sigma$ is simply the polynomial ring 
$
\C[z_\rho\,\big|\, \rho\in \Sigma(1)]
$ on $C^{\Sigma(1)}$,
which was introduced by D. Cox (and also called the {\it homogeneous coordinate ring}) in \cite{Cox95}.
Following \cite[\S 5.4]{CLS11}, the $\Phi$-homogenization of $p$ (as well as that of $\widehat p$) is the polynomial $\pp$ given by
 \begin{equation}\label{eq: Phihomogenizationofp}
\pp:=\sum_{m\in \Phi\cap M }c_m\cdot \prod_{\rho\in\Sigma(1)}z_\rho^{\langle m,u_\rho\rangle+
a_\rho}\in \C[z_\rho\,\big|\, \rho\in \Sigma(1)]
,
 \end{equation}
where $u_\rho$ and $a_\rho$ are as in \eqref{eq: dualfacetpresentation} (cf. also 
\eqref{eq: supporthyperplane}).  (More generally, for each 
$m\in \Phi\cap M$, the $\Phi$-homogenization of $\chi^m$ is the monomial $\displaystyle 
\prod_{\rho\in\Sigma(1)}z_\rho^{\langle m,u_\rho\rangle+
a_\rho}$.)
We remark that the inclusion 
in \eqref{eq: Phihomogenizationofp} follows from the inequalities in \eqref{eq: dualfacetpresentation}.
Upon polarizing  $\pp$ and similar to \eqref{eq: Pdef2}, we obtain the polynomial $\widetilde P
\in \C[z_\rho, \overline{w_\rho}\,\big|\, \rho\in \Sigma(1)]$ given by
 \begin{equation}\label{eq: PPhihomogenizationofp}
\PP(z,\overline w):=\sum_{m\in \Phi\cap M }c_m\cdot \prod_{\rho\in\Sigma(1)}
z_\rho^{\langle m,u_\rho\rangle+a_\rho}
\overline{
w_\rho^{\langle m,u_\rho\rangle+a_\rho}}
 \end{equation}
 for $z=(z_\rho)_{\rho\in\Sigma(1)},~w=(w_\rho)_{\rho\in\Sigma(1)}\in \C^{\Sigma(1)}$.
Consider the torus group $(\C^*)^{\Sigma(1)}$ acting on $\C^{\Sigma(1)}$ via coordinatewise multiplication, and let $G\subset (\C^*)^{\Sigma(1)}$ be the subgroup given by
\begin{align}
\label{eq:  defG}
G:&=
\{(t_\rho)_{\rho\in\Sigma(1)}\in
 (\C^*)^{\Sigma(1)}\,\big|\,\prod_{\rho\in \Sigma(1)}t_\rho^{\langle m,u_\rho\rangle}=1\text{ for all }m\in M\}.\\
\nonumber
&\cong \text{Hom}_{\Z}(\text{Cl}(X_\Sigma),\C^*),
\end{align}
where $\text{Cl}(X_\Sigma)$ denotes the Weil divisor group on $X_\Sigma$
(cf. e.g. \cite[\S 5.1]{CLS11} for more detailed discussion on the above isomorphism).  
Then it is known that $X_\Sigma $ is also given by
\begin{align}
\label{eq:  XSigmaasquotient}
X_\Sigma&=(\C^{\Sigma(1)}\setminus Z(\Sigma) )/G,\quad\text{where}\\
\label{eq:  irrelevantset}
Z(\Sigma):&=
\Big\{(z_\rho)_{\rho\in\Sigma(1)}\in
\C^{\Sigma(1)}\,\big|\,\prod_{\rho\in \Sigma(1)\setminus\sigma(1)}z_\rho=0\text{ for all }\sigma\in\Sigma(n)\Big\}
\end{align}
(cf. e.g. \cite[\S 5.1]{CLS11}).  We will denote the projection map associated to \eqref{eq:  XSigmaasquotient}
by 
\begin{equation}\label{eq: defXi}
\Xi: \C^{\Sigma(1)}\setminus Z(\Sigma)\to X_\Sigma.
\end{equation}
In light of \eqref{eq:  defL}, we let
$\chi^L:G\to \C^*$ be the group character given by
 \begin{equation}\label{eq: Gcharacter}
\chi^L((t_\rho)_{\rho\in\Sigma(1)})=
\prod_{\rho\in \Sigma(1)}t_\rho^{a_\rho}
\quad\text{for all }(t_\rho)_{\rho\in\Sigma(1)}\in G.
 \end{equation}
Together with \eqref{eq: Phihomogenizationofp}
and \eqref{eq:  defG}, one sees that
 \begin{equation}\label{eq: pfunctionaleqn}
\pp(g\cdot z)=\chi^L(g)\cdot \pp(z)\quad\text{for all }g\in G
\text{ and }z\in \C^{\Sigma(1)}
 \end{equation}
 (cf. e.g.  \cite[\S 5.2]{CLS11}).
 \begin{remark}\label{rem: triviallb}
One easily checks that $\Xi^*L$ is the trivial line bundle $\mathcal O_{ \C^{\Sigma(1)}\setminus Z(\Sigma)}$ on $\C^{\Sigma(1)}\setminus Z(\Sigma)$, and 
the $\Phi$-homogenization of a section $s\in  H^0(X_\Sigma, L)$ (arising from those of the $\chi^m$'s) is 
the holomorphic function $\widetilde s=\Xi^*s$ on $ \C^{\Sigma(1)}\setminus Z(\Sigma)$ which satisfies
the functional equation in \eqref{eq: pfunctionaleqn} (with $\pp$ replaced by $\widetilde{s}$).   
\end{remark}

\medskip Since $\Phi$ is an $n$-dimensional smooth lattice polytope, 
$X_\Sigma$ is covered by affine coordinate open subsets $\{U_\sigma\,\big|\, \sigma\in \Sigma(n)\}$, where 
\begin{equation}\label{eq:  Usigma}
U_\sigma:=\C^{\sigma(1)}
=\{(z_\rho)_{\rho\in \sigma(1)}\,\big|\, z_\rho\in\C\}
\cong \C^n\quad\text{for each }\sigma\in \Sigma(n),
\end{equation}
and one has a natural holomorphic map $\phi_\sigma:U_\sigma\to \C^{\Sigma(1)}$ given by 
\begin{equation}\label{eq: phisigma}
(\phi_\sigma(z))_{\rho}:=\begin{cases}z_\rho\quad\text{if }\rho\in \sigma(1),
\\
1\quad\text{if }\rho\in\Sigma(1)\setminus \sigma(1)
\end{cases}
\end{equation}
for $z=(z_\rho)_{\rho\in \sigma(1)}$
(cf. e.g.  \cite[\S 5.2]{CLS11}).
It is easy to see that for each $\sigma\in \Sigma(n)$, one obtains a holomorphic trivialization  
\begin{equation}\label{eq: LtrivialUalpha}
L\big|_{U_\sigma}\cong
U_\sigma\times \C\quad\text{given by}\quad
s(z)\longleftrightarrow(z,\widetilde{s}(\phi_\sigma(z)))
\end{equation}
for $z\in U_\sigma\text{ and }s\in H^0(X_\Sigma, L)$.

\medskip 
 Our main result in this section is the following proposition:
\begin{proposition} \label{prop: P3properties}
Let $p \in \R[x_1^{\pm 1},\dots,x_n^{\pm 1}]$ be a Laurent polynomial whose Newton polytope $\Phi$ is an
$n$-dimensional smooth lattice polytope. Let $\pp$ be the $\Phi$-homogenization of $p$, and let 
$\widehat P\in  H^0(X_\Sigma\times \overline{X_\Sigma}, \pi_1^*L\otimes \pi_2^*\overline{L})$ be the associated Hermitian algebraic function on $L$ as given in \eqref{eq: Pdef2}.
If $\pp$
    satisfies (Pos1), (Pos2), and (Pos3), then the following statements hold:
    \begin{enumerate}[(i)]
\item \label{prop: Pposdef}  $\widehat P$ is positive (as a Hermitian algebraic function on $L$).
\item \label{prop: positivecur}  The curvature $(1,1)$-form $\Theta_{h_{\widehat P}}$ on $X_\Sigma$ is negative definite.   
\item \label{prop: PSGCS} $\widehat P$ satisfies the SGCS inequality.
\end{enumerate}
\end{proposition}

Throughout the remainder of this section, which is devoted to the proof of the above proposition, we let 
$p \in \R[x_1^{\pm 1},\dots,x_n^{\pm 1}]$ be a Laurent polynomial whose Newton polytope $\Phi$ is an
$n$-dimensional smooth lattice polytope, $\pp$ be the $\Phi$-homogenization of $p$, and 
$\widehat P\in  H^0(X_\Sigma\times \overline{X_\Sigma}, \pi_1^*L\otimes \pi_2^*\overline{L})$ be the associated Hermitian algebraic function on $L$ as above.
First we have

\begin{proposition}\label{prop: Pos1Pos3impliesPpd} 
If $\pp$ satisfies (Pos1) and (Pos3), 
    then 
\begin{enumerate}[(i)]
\item \label{prop: ppositive} $\pp(x)>0 $ for all $x\in\R_+^{\Sigma(1)} \setminus (\R_+^{\Sigma(1)}\cap Z(\Sigma))$, and 
\item \label{prop: Ppdf} $\widehat P$ is positive (as a Hermitian algebraic function on $L$).
\end{enumerate}
\end{proposition}

\begin{proof} 
To prove (i), we let $x=(x_\rho)_{\rho\in\Sigma(1)} \in \R_+^{\Sigma(1)} \setminus (\R_+^{\Sigma(1)}\cap Z(\Sigma))$ be given, and let $\pp$ be as in \eqref{eq: Phihomogenizationofp}.
We consider the following two cases:
\par\noindent
$\underline{Case~(a)}$:~$ (G\cap U(1)^{\Sigma(1)} )\cdot x 
\neq
U(1)^{\Sigma(1)}\cdot x$.  In this case, we take $z\in 
(U(1)^{\Sigma(1)}\cdot x)\setminus ((G\cap U(1)^{\Sigma(1)} )\cdot x )$ (and thus
$(|z_\rho|)_{\rho\in\Sigma(1)}=x$), then it follows from (Pos3)
that $\pp(x)>|\pp(z)|\geq 0$.   
\par\noindent
$\underline{Case~(b)}:~$ In this remaining case, we have
\begin{equation}\label{eq: GRorbitx} 
(G\cap U(1)^{\Sigma(1)} )\cdot x =U(1)^{\Sigma(1)}\cdot x.
\end{equation} 
Consider the holomorphic map $\phi: (\C^*)^{\Sigma(1)}\to  \C^{\Sigma(1)}$ given by 
$\phi(t)=t\cdot x$ for $t\in (\C^*)^{\Sigma(1)}$, where $t\cdot x  $ denotes coordinatewise multiplication.  Let $\Sigma^\prime(1):=\{\rho\in\Sigma(1)\,\big|\, x_\rho\neq 0\}$.  It is easy to see that $\phi ((\C^*)^{\Sigma(1)})=(\C^*)^{\Sigma^\prime(1)}$, where $(\C^*)^{\Sigma^\prime(1)}$ is 
regarded as a submanifold in $(\C^*)^{\Sigma(1)}$ in the obvious manner.  In particular, one has
\begin{equation}\label{eq: partialphi}
\partial\phi(T^{1,0}_{  (1,1,\cdots,1)}(\C^*)^{\Sigma(1)})=T^{1,0}_{ x}     (\C^*)^{\Sigma^\prime(1)}\cong\C^{\Sigma^\prime(1)}.
\end{equation}
It is easy to see that the complex linear span of  $T_{x}((\C^*)^{\Sigma(1)})$  is simply 
$\C^{\Sigma^\prime(1)}$ under the isomorphism in \eqref{eq: partialphi}.  Together with 
\eqref{eq: GRorbitx}, it follows that 
\begin{equation}\label{eq: partialphiG}
\partial\phi(T^{1,0}_{  (1,1,\cdots,1)}G)=\C^{\Sigma^\prime(1)}.
\end{equation} 
Note that $\Xi(G\cdot x)=\Xi(x)$ is a single point in $X_\Sigma$
(cf. \eqref{eq: defXi}).  Together with \eqref{eq: GRorbitx}  and \eqref{eq: partialphiG}, it follows that
$\Xi((\C^*)^{\Sigma(1)}\cdot x)$ is also a connected $0$-dimensional analytic set in $X_\Sigma$, and thus 
\begin{equation}\label{eq: singleorbit}\Xi((\C^*)^{\Sigma(1)}\cdot x)=\Xi(x), \text{  implying that }T_N\cdot \Xi(x)=\Xi(x)\text{ in }X_\Sigma.
\end{equation}
Note that $x=(x_\rho)_{\rho\in\Sigma(1)}\in \C^{\Sigma(1)}\setminus Z(\Sigma)$, and thus there exists $\sigma\in \Sigma(n)$ such that $x_\rho\neq 0$ for all $\rho\in \Sigma(1)\setminus \sigma(1)$ and $\Xi (x)\in U_\sigma$ (cf. \eqref{eq:  irrelevantset} and \eqref{eq:  Usigma}).  From the second equality in 
\eqref{eq: singleorbit}, one easily sees that 
\begin{equation}\label{eq: Xi0}
\Xi(x)=0 \text{  in }U_\sigma\cong \C^n \implies x_\rho=0\quad\text{for all }\rho\in\sigma(1).
\end{equation}
It follows readily  that the set of points $y$ in $\R_+^{\Sigma(1)} \setminus (\R_+^{\Sigma(1)}\cap Z(\Sigma))$ that correspond to Case (a) (so that $\pp(y)>0$) form a dense open subset in $\R_+^{\Sigma(1)} \setminus (\R_+^{\Sigma(1)}\cap Z(\Sigma))$.  Together with the continuity of $\pp$
on $\R_+^{\Sigma(1)}$, it follows that we have 
\begin{equation}\label{eq: ppgeq0}
\pp(x)\geq 0.
\end{equation}
Furthermore, let $t=(t_\rho)_{\rho \in \Sigma(1)}\in (\C^*)^{\Sigma(1)}$ be given by 
$t_\rho=1/x_\rho$ if $\rho\in \Sigma(1)\setminus\sigma(1)$ and $t_\rho=1 $ if $\rho\in \sigma(1)$.
Together with \eqref{eq: Xi0}, one sees that 
\begin{equation}\label{eq: tdotx0}
t\cdot x=e^{(\sigma)},
\end{equation}
where $e^{(\sigma)}$ is as in \eqref{eq: defesigma}.  From \eqref{eq: tdotx0} and the first equality in 
\eqref{eq: singleorbit}, one sees that $\Xi( e^{(\sigma)} )=  \Xi(x)$, and thus there exists $g\in G$ such that
$g\cdot x=e^{(\sigma)}$.  Then from \eqref{eq: pfunctionaleqn}, one has
\begin{equation}\label{eq: pfunctionaleqn2}
\pp(e^{(\sigma)})=\pp(g\cdot x)=\chi^L(g)\cdot \pp(x).
 \end{equation}
Together with (Pos1) and the fact that $\chi^L(g)\in \C^*$, one sees that that $\pp(x) \neq 0$, which 
together with \eqref{eq: ppgeq0}, imply that $\pp(x) > 0$. Thus we have finished the proof of (i).
We proceed to prove (ii).  For each $\sigma\in \Sigma(n)$, we let $U_\sigma$ be as in 
\eqref{eq:  Usigma}.  To prove that the Hermitian algebraic function $\widehat{P}$ on $L$ is positive, it suffices to establish the positivity of $\widehat{P}$ on each such $U_\sigma$.  
In terms of the trivialization of $L\big|_{U_\sigma}$ given in \eqref{eq: LtrivialUalpha} and $z=(z_\rho)_{\rho\in\sigma(1)}\in U_\sigma
\cong \C^n$, one sees that $\widehat{P}(z,\overline{z})$ is given by $\PP(\phi_\sigma(z),
\overline{\phi_\sigma(z)})$, where $\phi_\sigma$ is as in \eqref{eq: phisigma}.
Together with \eqref{eq: Phihomogenizationofp} and \eqref{eq: PPhihomogenizationofp}, one sees that
\begin{align}\label{eq: PPequalsPPx1}
\PP(\phi_\sigma(z),
\overline{\phi_\sigma(z)})&=\pp(x),\quad\text{where }x=(x_\rho)_{\rho\in\Sigma(1)}\text{ is given by}\\
\label{eq: PPequalsPPx2}
x_\rho&=\begin{cases}|z_\rho|^2\quad\text{if }\rho\in\sigma(1),
\\
1\quad\text{if }\rho\in\Sigma(1)\setminus\sigma(1).
\end{cases}
 \end{align}
Since $x_\rho\neq 0$ for all $\rho\in\Sigma(1)\setminus\sigma(1)$, it follows that $x\notin Z(\Sigma)$ (cf. 
\eqref{eq:  irrelevantset}).  Thus, 
$x\in \R_+^{\Sigma(1)} \setminus (\R_+^{\Sigma(1)}\cap Z(\Sigma))$, which together with \eqref{prop: ppositive},  imply that $\pp(x)>0$.  This finishes the proof of Proposition \ref{prop: Pos1Pos3impliesPpd}(ii).
\end{proof}

Next we recall a result of De Angelis \cite{Deangelis94}.  For $\ell \geq 1$, we denote the interior of 
$\R^\ell_+$ by $(\R^\ell_+)^\circ:=\{(s_1,\dots, s_\ell)\in \R^\ell\,\big|\, s_i>0,~i=1,\dots,\ell\}$.  
Let $f(s)  
\in \R[s_1, \ldots, s_\ell]$ be 
	such that $f(s) > 0$ for all $s =(s_1,\dots, z_\ell)\in (\R^\ell_+)^\circ$.  
We associate to $f$ the $\ell \times \ell$ matrix-valued function $J_f
: (\R^\ell_+)^\circ\to \R^{\ell^2}$ whose components are given by
\begin{align} \label{eq: defineJ}
J_f(s)_{ij} :&= s_j \cdot \frac{\partial}{\partial s_j}\Big(s_i\cdot\frac{\partial }{\partial s_i}\big(\log f\big)\Big)(s)
\\
&=s_is_j\frac{\partial^2 } {\partial s_i\partial s_j  }\big(\log f\big)(s)+\delta_{ij}\cdot s_j\cdot
\frac{\partial}{\partial s_i}\big(\log f\big)(s)
\notag
\end{align}
for $s\in (\R^{\ell}_+)^\circ,~1\leq i,j\leq \ell$.  Here $\delta_{ij}$ denotes the Kronecker delta. i.e., $\delta_{ij}=1 $ (resp.
$0$) if $i=j$ (resp. $i\neq j$).  Next we introduce a change of variables, and consider 
the function $f^\sharp:\R^{\ell}\to \R$ associated to $f$ given by 
\begin{equation}\label{eq: tildef}
f^\sharp(t)=f(e^{t_1}, \ldots, e^{t_\ell}) \quad \text{for } t = (t_1,\ldots, t_\ell) \in \R^{\ell}.
\end{equation}
Using \eqref{eq: defineJ}, one easily checks that the Hessian matrix of 
$\log f^\sharp$ coincides with $J_f$, i.e., one has
\begin{equation}\label{eq: HessJlogf}
\frac{\partial^2 } {\partial t_i\partial t_j  }\big(\log f^\sharp\big)(t) =J_f(e^{t_1}, \ldots, e^{t_\ell})_{ij}
\end{equation}
for all $t=(t_1,\dots, t_\ell)\in\R^{\ell},~1\leq i,j\leq \ell$.   We recall the following result:
	
\begin{lemma}[De Angelis {\cite[Theorem 6.11]{Deangelis94}}] \label{lem: positiveHessianInInterior}
Let $f(s) \in \R[s_1, \ldots, s_\ell]$ be such that $f(s) > 0$ for all $s \in (\R^{\ell}_+)^\circ$.
Suppose that there exists an open neighborhood $V$ of $(\R^{\ell}_+)^\circ$ in $(\C ^*)^\ell$
		such that $|f(z)| \le f(|z_1|, \ldots, |z_\ell|)$ for all $z = (z_1, \ldots, z_\ell) \in V$, and the 
		Newton polytope $\Phi(f)$ of $f$ has affine dimension $\ell$.
Then the $\ell \times \ell$ matrix $J_f(s)$ is positive definite for all $s \in (\R^{\ell}_+)^\circ$.
\end{lemma}

We let $p, \,\widehat{p}, \,\widehat{P},\, \pp,\, \PP $ be as before. Let 
$\mathfrak{S}_n$ denote the group of
permutations of the coordinate functions on $\C^n$.  For each $1\leq \ell\leq n$, each 
$\sigma\in\Sigma(n)$ and each $\tau\in \mathfrak{S}_n$, we 
let $\pp_{\ell,\sigma,\tau}\in \R[s_1,\dots, s_\ell]$ be
given by
\begin{equation}\label{eq: definepmsigma}        
\pp_{\ell,\sigma,\tau}:=\pp(\phi_\sigma(\tau(s_1,\cdots, s_\ell,0,\dots,0))),
\end{equation} 
where $\phi_\sigma$ is as in \eqref{eq: phisigma}.

\begin{lemma} \label{lem: Zspan}  
Let $p \in \R[x_1^{\pm 1},\dots,x_n^{\pm 1}]$ be such that its Newton polytope $\Phi$ is an
$n$-dimensional smooth lattice polytope.
\par\noindent
(i)  If $\pp$ satisfies (Pos1), then for each $1\leq \ell\leq n$, each 
$\sigma\in\Sigma(n)$ and each $\tau\in \mathfrak{S}_n$, the Newton polytope 
$\Phi(\pp_{\ell,\sigma,\tau})$
of $\pp_{\ell,\sigma,\tau}$ has affine dimension $\ell$. 
\par\noindent
(ii)  If $\pp$ satisfies (Pos1) and (Pos2), then for each $1\leq \ell\leq n$, each 
$\sigma\in\Sigma(n)$ and each $\tau\in \mathfrak{S}_n$, the set
$
S_{\pp_{\ell,\sigma,\tau}} : =\{m - m^\prime \, \big|\,m,\,m^\prime \in \Log (\pp_{\ell,\sigma,\tau})\}
$
generates $\Z^{\ell}$ as a 
$\Z$-module.
\end{lemma}   

\begin{proof}  As the proofs of the lemma for all the $\pp_{\ell,\sigma,\tau}$'s are the same, 
 we will only prove the lemma for the case when $\tau$ is the identity permutation, so that
 $
\pp_{\ell,\sigma,\tau}(s_1,\dots,s_\ell)=\pp(\phi_\sigma(s_1,\cdots, s_\ell,0,\cdots,0))$.  
Fix a cone $\sigma\in\Sigma(n)$ and an integer $\ell$ satisfying $1\leq\ell\leq n$.  Denote the cardinality of $\Sigma(1)$ by $|\Sigma(1)|$.  For notational convenience, we write 
\begin{equation}\label{eq: sigmaSigma1}
\sigma(1)=\{\rho_i\,\big|\, i=1,\cdots,n\}\text{  and  }\Sigma(1)\setminus\sigma(1)=\{\rho_i
\,\big|\, i=\ell+1,\cdots,|\Sigma(1)|\}
\end{equation}
such that $\phi_{\sigma}$ is given by 
\begin{equation}\label{eq: phisigma1n}
\phi_{\sigma}(z_{\rho_1},\cdots,z_{\rho_n})=(z_{\rho_1},\cdots,z_{\rho_n},1,\cdots, 1)
\end{equation}
(cf. \eqref{eq: phisigma}).  
For each $1\leq j\leq n$, we let $H_j\subset M_\R$ be the hyperplane given by 
\begin{equation}\label{eq:defHj}
H_j:=\{m\in M_\R\,\big|\, \langle m, u_{\rho_j}\rangle=-a_{\rho_j}\}.
\end{equation}
Then from \eqref{eq: Phihomogenizationofp}, one easily checks that
\begin{equation}\label{eq: ppellsigma1ell}
\pp_{\ell,\sigma,\tau}(s_1,\dots,s_\ell)=\sum_{m\in \Phi\cap M \cap H_{\ell+1}\cap\cdots\cap H_n}c_m\cdot \prod_{1\leq i\leq \ell}s_i^{\langle m,u_{\rho_i}\rangle+
a_{u_{\rho_i}}}
.
\end{equation}
Let $U_n $ (resp. $U_\ell$) be the $n\times n$ ( resp. $\ell\times n$) matrix such that its $i$-th row is given by $u_{\rho_i}$ (written as a row matrix), and let $A_n $ (resp. $A_\ell$) be the column vector with $n$ (resp. $\ell$) rows such that its $i$-th entry is $a_{\rho_i}$.   Let $\eta_n: M_\R\to \R^n$ and $\eta_\ell:M_\R\to \R^\ell$ be given by 
\begin{equation}\label{eq: defetan}
\eta_n(m)=U_n m+A_n \text{  and  }\eta_\ell(m)=U_\ell m+A_\ell \quad\text{for }m\in M_\R.
\end{equation}
From \eqref{eq: ppellsigma1ell}, one sees that 
\begin{equation}\label{eq: ppnewton}
\Phi(\pp_{\ell,\sigma,\tau})=\text{ convex hull of }\{\eta_\ell(m)\,\big|\, m\in \Phi\cap M \cap H_{\ell+1}\cap\cdots\cap H_n\}.
\end{equation}  
From \eqref{eq:defHj} and \eqref{eq: defetan}, one easily sees  that for each $m\in \Phi\cap M \cap H_{\ell+1}\cap\cdots\cap H_n$, one has 
\begin{equation}\label{eq: compareetanell}
(\eta_n(m))_j=\begin{cases}(\eta_\ell(m))_j\quad\text{if }1\leq j\leq\ell,\\
0\quad\text{if }\ell+1\leq j\leq n,
\end{cases}
\end{equation} 
where $(\eta_n(m))_j$ (resp. $(\eta_\ell(m))_j$) denotes the $j$-th entry of $\eta_n(m)$ (resp. $\eta_\ell (m)$).
Since $\pp$ satisfies (Pos1), one sees from \eqref{eq: defesigma}, 
\eqref{eq: Phihomogenizationofp} and \eqref{eq: ppellsigma1ell} that 
\begin{equation}\label{eq: 0inNewtonpp}
\pp_{\ell,\sigma,\tau}(0,\cdots,0)=\pp(e^{(\sigma)})>0,\quad\text{so that  }(0,\cdots,0)\in \Phi(\pp_{\ell,\sigma,\tau}).
\end{equation}
Let $v\in\Phi\cap M$ be the vertex of $\Phi$ that corresponds to $\sigma$ under the Orbit-Cone Correspondence.   Then it is easy to see that $\{v\}=\cap_{j=1}^n H_j$, so that
$\eta_n(v)=0$ and thus $\eta_\ell(v)=0$.  
Since $\Phi$ is a smooth polytope, it follows that $\sigma$ is a smooth cone, which implies that the matrices $U_n$ and $U_\ell$ are of rank $n$ and $\ell$ respectively; furthermore, the $\R$-span of the set $\{m-v\,\big|\, m\in \Phi\cap M \cap H_{\ell+1}\cap\cdots\cap H_n\}$ has dimension $\ell$ (see e.g. \cite[\S 2.4]{CLS11}).  Since $U_n$ is of rank $n$, it follows that  the $\R$-span of the set $\{\eta_n(m)\,\big|\, m\in \Phi\cap M \cap H_{\ell+1}\cap\cdots\cap H_n\}$ also has dimension $\ell$ (noting that $\eta_n(m)=\eta_n(m-v)$).  Together with 
\eqref{eq: compareetanell}, one sees  that the $\R$-span of the set $\{\eta_\ell(m)\,\big|\, m\in \Phi\cap M \cap H_{\ell+1}\cap\cdots\cap H_n\}$ also has dimension $\ell$.  Combining this with \eqref{eq: ppnewton} and 
\eqref{eq: 0inNewtonpp}, one sees that $\Phi(\pp_{\ell,\sigma,\tau})$ has affine dimension $\ell$, and this finishes the proof of (i).  
 We proceed to prove (ii) and adopt the same notation as above.  Write $\pp\in\R[x_{\rho_i}\,\big|\, i=1,\cdots, |\Sigma(1)|\}$ (cf. \eqref{eq: sigmaSigma1}).
For each 
$1\leq i\leq \ell$, one 
easily checks from \eqref{eq: defesigma}, 
\eqref{eq: Phihomogenizationofp} and \eqref{eq: ppellsigma1ell} that
\begin{equation}\label{eq: partialpmsigma}
\frac{\partial \pp_{\ell,\sigma,\tau}}{\partial s_i}(0,\dots,0)=\frac{\partial \pp}{\partial x_{\rho_i}}(e^{(\sigma)})>0,
\end{equation}
where the inequality holds since $\pp$ satisfies (Pos2) and $e^{(\sigma)}\in F_{\rho_i}(\R_+^{\Sigma(1)}) \setminus
  											 (Z(\Sigma)\cap F_{\rho_i}(\R_+^{\Sigma(1)}) )$ (note that $e^{(\sigma)}\notin 
  			Z(\Sigma)$ since $e^{(\sigma)	}_\rho=1\neq 0$ for
  			each $\rho\in \Sigma(1)\setminus\sigma(1)$ (cf.						
  											 \eqref{eq:  irrelevantset})). 
This implies that $\Log (\pp_{\ell,\sigma,\tau})(\subset \Z^{\ell})$ contains the points $(1,0,\dots,0)$, $\ldots$, $(0,\dots, 0,1)$, and so does $S_{\pp_{\ell,\sigma,\tau}}$ 
(since $\Log (\pp_{\ell,\sigma,\tau})$ also contains $(0,\cdots, 0)$ as shown in \eqref{eq: 0inNewtonpp}). 
It follows that $S_{\pp_{\ell,\sigma,\tau}}$ generates $\Z^{\ell}$ as a 
$\Z$-module.
\end{proof} 

\medskip
If $\pp$ satisfies (Pos1) and (Pos3), then it follows from Proposition
\ref{prop: Pos1Pos3impliesPpd} that $\widehat P$ is a positive Hermitian algebraic function on $L$, and thus as discussed in Section \ref{sec: CDThm}, 
$\widehat P$ induces a Hermitian metric $h_{\widehat{P}}$ on $L^*$ given by $h_{\widehat{P}}(v,\overline{w})=\widehat{P}(v,\overline{w})$ for $v,w\in L^*_x$, $x\in X_\Sigma$.   Here $
L^*_x$ denotes the fiber of $L^*$ at $x$.

\begin{proposition} \label{prop: StronglyPseudoconvex}
If $\pp$ satisfies (Pos1), (Pos2) and (Pos3),
    then the curvature $(1,1)$-form $\Theta_{h_{\widehat{P}}}$ on $X_\Sigma$ is negative definite.
\end{proposition}

\begin{proof} Take any cone $\sigma\in\Sigma(n)$, and let $U_\sigma\cong\C^n$ be the affine coordinate 
subset of $X_\Sigma$ 
with coordinate functions $z_\rho$, $\rho\in \sigma(1)$, as in \eqref{eq:  Usigma}.  Let
$\xi\in
H^0(U_\sigma, L\big|_{U_\sigma})$ be the non-vanishing section so that 
$\xi(z) \longleftrightarrow (z,1)$ with respect to  the holomorphic trivialization $L\big|_{U_\sigma}\cong U_\sigma\times \C$ in \eqref{eq: LtrivialUalpha}.  Let $\xi^*\in
H^0(U_\sigma, L^*\big|_{U_\sigma})$ be the dual non-vanishing section, so that $\langle \xi,\xi^*\rangle\equiv 1$ on $U_\sigma$.  Here $\langle \xi,\xi^*\rangle$ denotes the natural pointwise pairing
between $L$ and $L^*$.  
Upon following the notation in \eqref{eq: sigmaSigma1}, \eqref{eq: phisigma1n} and writing $\pp\in\R[x_{\rho_i}\,\big|\, i=1,\cdots, |\Sigma(1)|\}$ as in \eqref{eq: partialpmsigma}, 
one easily sees from \eqref{eq: Phihomogenizationofp}, \eqref{eq: PPhihomogenizationofp}, \eqref{eq: phisigma} \eqref{eq: LtrivialUalpha}  
that
\begin{equation}\label{eq: hxixi}
h_{\widehat{P}}(\xi^*,\overline{\xi^*})(z)=\PP(\phi_\sigma(z),\overline{\phi_\sigma(z)})=\pp(|z_{\rho_1}|^2,\cdots,
|z_{\rho_n}|^2,1,\cdots, 1)
\end{equation}
for $z=(z_{\rho_1},\cdots, z_{\rho_n})\in U_\sigma$.
Write 
$\displaystyle\Theta_{h_{\widehat P}}\big|_{U_\sigma}=\sum_{1\leq i,j\leq n}(\Theta^{(\sigma)}_{h_{\widehat P}})_
{i\overline{j}} dz_{\rho_i}\wedge d\overline{z}_{\rho_j}$.
Then from \eqref{eq: hxixi}, one has, for $z=(z_{\rho_1},\cdots, z_{\rho_n})\in U_\sigma$ and $1\leq i,j\leq n$, 
\begin{align}\label{eq: thetahPij}
&\qquad (\Theta^{(\sigma)}_{h_{\widehat P}}(z))_
{i\overline{j}}
\\
\nonumber
&=-\frac{\partial^2( \log h_{\widehat{P}}(\xi^*,\overline{\xi^*}))} {\partial z_{\rho_i}\partial\conj{z_{\rho_j}}}(z,\conj{z})
\\
\nonumber
&=-\Big(
z_{\rho_j}\conj{z}_{\rho_i}\cdot\frac{\partial^2( \log \pp)} {\partial x_{\rho_i}\partial x_{\rho_j}}
+\delta_{ij}\cdot\frac{\partial ( \log \pp)}{\partial x_{\rho_i}}\Big)
(|z_{\rho_1}|^2,\cdots,
|z_{\rho_n}|^2,1,\cdots, 1).
\end{align}
Take a point $z^*=(z_{\rho_1}^*,\dots, z_{\rho_n}^*)\in U_\sigma $, and let
$\ell$ be the number of non-zero $z_{\rho_i}^*$'s for $1\leq i\leq n$ (so that $0\leq \ell\leq n$).
By permuting the $z_{\rho_i}$'s, we will assume without loss of generality that
$z_{\rho_i}^*\neq 0$ for each $1\leq i\leq \ell$ and $z_{\rho_{\ell+1}}^*=\cdots=z_{\rho_{n}}^*=0$. 
Now we take a tangent vector $0\neq v=
v_1\dfrac{\partial}{\partial z_{\rho_1}}+\cdots +v_n\dfrac{\partial}{\partial  z_{\rho_n}}\in T_{z^*}U_\sigma$.
By using \eqref{eq: thetahPij} and \eqref{eq: defineJ} (with $f=\pp_{\ell,\sigma,\text{Id}}$ where 
$\text{Id}$ denotes the identity
permutation (cf. \eqref{eq: definepmsigma}), and 
$s=(|z_{\rho_1}^*|^2,\dots,|z_{\rho_\ell}^*|^2)$), one easily checks that 
\begin{align}\label{eq: twoHessian}
\sum_{1\leq i,j\leq n}&\conj{v_j}\cdot
(\Theta^{(\sigma)}_{h_{\widehat P}}(z^*))_
{i\overline{j}}\cdot v_i =-A_1-A_2,\quad\text{where}\\
\notag
A_1:&=\sum_{1\leq i,j\leq \ell} \dfrac{\conj{ v_j}}{ \conj{z_{\rho_i}^*} }
 \cdot J_{\pp_{\ell,\sigma,\text{Id}}}(|z_{\rho_1}^*|^2,\dots,|z_{\rho_\ell}^*|^2)_{ij}\cdot \dfrac{ v_i }{ {z_{\rho_i}^*} }
\quad\text{and}
  \\
  \notag
A_2:= & \sum_{\ell+1\leq i\leq n}|v_i|^2\cdot\dfrac{\partial (\log \pp)}{\partial x_{\rho_i}}
 (|z_{\rho_1}^*|^2,\dots,|z_{\rho_{\ell}}^*|^2,0,\cdots,0,1,\cdots,1).
\end{align}
Here $A_1$ (resp. $A_2$) is taken to be zero if $\ell=0$ (resp. $\ell=n$).
Note that $  (|z_{\rho_1}^*|^2,\dots,|z_{\rho_{\ell}}^*|^2,0,\cdots,0,1,\cdots,1)\in F_{\rho_i}(\R_+^{\Sigma(1)})$ for each $\ell+1\leq i\leq n$.  
Hence from (Pos2), 
we see that $A_2\geq 0$, and
\begin{equation}\label{eq: A2greaater0}
A_2>0\quad\text{whenever }\ell<n\text{ and }(v_{\ell+1},\dots, v_{n})\neq (0,\dots,0).
\end{equation}
From  (Pos3) (for the set $ \C^{\Sigma(1)} \setminus (Z(\Sigma)\cup (G\cap U(1)^{\Sigma(1)})\cdot \R_+^{\Sigma(1)})$, which is easily seen to be dense in $\C^{\Sigma(1)}$) and the 
continuity of $\pp$, one easily sees that $|\pp(z) |\leq \pp(|z_{\rho_1}|,\dots,|z_{\rho_{ |\Sigma(1)|}}|)$ for all $z=(z_{\rho_1},\dots,
z_{\rho_{ |\Sigma(1)|}}) \in 
\C^{\Sigma(1)}$, which implies readily that
\begin{equation}\label{eq:  pinequalitycn}
|\pp_{\ell,\sigma,\text{Id}}(z)|\leq \pp_{\ell,\sigma,\text{Id}}(|z_1|,\cdots, |z_\ell|)\quad
\text{for all }z=(z_1,\cdots, z_\ell)\in \C^\ell.
\end{equation}
Together with Lemma \ref{lem: Zspan}, it follows that one can apply Lemma \ref{lem: positiveHessianInInterior} (with $f=\pp_{\ell,\sigma,\text{Id}}$)
to conclude that $A_1\geq 0$, and
\begin{equation}\label{eq: A1greaater0}
A_1>0\quad\text{whenever }\ell>0\text{ and }(v_1,\dots, v_{\ell})\neq (0,\dots, 0).
\end{equation}
Together with \eqref{eq: A2greaater0},
one easily concludes that $A_1+A_2>0$ 
in each of the three cases when $\ell=0$, $1\leq \ell<n$ or $\ell=n$.  
Upon varying $z^*\in U_\sigma$ and using the fact that $\{U_\sigma\,\big|\, \sigma\in\Sigma(n)\}$ covers $X_\Sigma$, one
concludes that $\Theta_{h_{\widehat P}}$ is positive definite on $X_\Sigma$.
\end{proof}

\begin{lemma}\label{lem: gx}
 Let  
 $\sigma\in\Sigma(n)$ and $x=(x_{\rho})_{\rho\in\Sigma(1)|} \in \R_+^{\Sigma(1)} $ be such that $x_\rho>0$ for all $\rho\in \Sigma(1)\setminus\sigma(1)$.  Then 
there exist $g\in G\cap  (\R_+^{\Sigma(1)})^\circ $ and $s\in \R_+^{\sigma(1)})$ such that
$g\cdot x=\phi_\sigma(s)$, where $\phi_\sigma$ is as in \eqref{eq: phisigma}.  Furthermore, if $x\in 
(\R_+^{\Sigma(1)})^\circ$, then $s\in (\R_+^{\sigma(1)})^\circ$.
\end{lemma}

\begin{proof} We write $\sigma(1)$ and $\Sigma(1)$ as in \eqref{eq: sigmaSigma1}.
Let  
$x=(x_{\rho_i})_{1\leq i\leq |\Sigma(1)|} \in \R_+^{\Sigma(1)}  $ be such that $x_{\rho_i}>0$ for each $n+1\leq i\leq |\Sigma(1)|$.  We need to find $g=(t_{\rho_i})_{1\leq i\leq |\Sigma(1)|}\in G$ with
each $t_{\rho_i}>0$ and $s=(s_{\rho_i})_{1\leq i\leq n}\in \R_+^{\sigma(1)}$ such that
\begin{equation}\label{eq: solveeqngx}
t_{\rho_i}\cdot x_{\rho_i}=\begin{cases}s_{\rho_i}\quad\text{if }1\leq i\leq n,\\
1\quad\text{if }n+1\leq i\leq |\Sigma(1)|.
\end{cases}
\end{equation}
For this purpose, we first let $t_{\rho_i}:=1/x_{\rho_i}>0$ for each $n+1\leq i\leq |\Sigma(1)|$, and let
$T$ be the $( |\Sigma(1)|-n)$-column vector whose $i$-th entry is $\log t_{\rho_{n+i}}$.  Then we let $A$ 
(resp. $B$) be the $n\times n$ (resp. $n\times ( |\Sigma(1)|-n)$) matrix whose $i$-th column is $u_{\rho_i}$ (resp. $u_{\rho_{n+i}}$), where the $u_{\rho_i}$'s are as in 
\eqref{eq: dualfacetpresentation}.  Since $\sigma$ is an $n$-dimensional smooth cone, it follows that the matrix
$A$ is non-singular.  Now for each $1\leq i\leq n$, we let $t_{\rho_i}>0$ be the number such that
$\log t_{\rho_i}$ is the $i$-th entry of the column vector $-A^{-1}BT$, so that one has
\begin{equation}\label{eq: logtiu}
\sum_{i=1}^{|\Sigma(1)|}\log t_{\rho_i}\cdot u_{\rho_i}=0 \quad \text{in }\R^n.
\end{equation}
Then upon taking inner product of both sides of \eqref{eq: logtiu} with each $m\in M$ and exponentiating the resulting expressions, one easily sees 
from \eqref{eq:  defG} that $g=(t_{\rho_i})_{1\leq i\leq  |\Sigma(1)}\in G$ with each $t_{\rho_i}>0$.  Finally we let $s_{\rho_i}:=t_{\rho_i}\cdot x_{\rho_i}$ for each $1\leq i\leq n$.  Then one easily sees that  \eqref{eq: solveeqngx} is satisfied. Finally the last statement of Lemma \ref{lem: gx} is obvious.
\end{proof}

 \begin{remark}\label{rem: gx111}
 Let  
 $\sigma\in\Sigma(n)$ and $z=(z_{\rho})_{\rho\in\Sigma(1)|} \in \C^{\Sigma(1)} $ be such that $z_\rho\neq 0$ for all $\rho\in \Sigma(1)\setminus\sigma(1)$.  Then by following the proof of Lemma 
 \ref{lem: gx}, one can easily show that 
there exist $g\in G $ and $s\in \C^{\sigma(1)}$ such that
$g\cdot z=\phi_\sigma(s)$, where $\phi_\sigma$ is as in \eqref{eq: phisigma}. 
 \end{remark}

\begin{lemma}\label{lem: pconvex}
Suppose $\pp$ satisfies (Pos1), (Pos2) and (Pos3).  
 Let  
$k,\ell\in\N$ be such that $1\leq k\leq\ell\leq n$, $\sigma\in\Sigma(n)$, $\tau\in \mathfrak{S}_n$, and $\pp_{\ell,\sigma,\tau}\in \R[s_1,\dots, s_\ell]$
be as in \eqref{eq: definepmsigma}.  
\par\noindent
(i)  Let $(s_1,\cdots, s_k)\in (\R_+^k)^\circ$, i.e., $s_i>0$ for all $1\leq i\leq k$.  Then one has
\begin{equation}\label{eq: plincreasing}
\pp_{\ell,\sigma,\tau}(s_1,\cdots,s_{k-1},0,\cdots,0)<\pp_{\ell,\sigma,\tau}(s_1,\cdots,s_{k-1},s_k,0,\cdots,0).
\end{equation}
\par\noindent
(ii)  Let $s = (s_i)_{1\leq i\leq \ell}\neq s^\prime = (s_i^\prime)_{1\leq i\leq \ell} \in (\R_+^{\ell})^\circ$.  Then one has
\begin{equation} \label{eqlem: convexEqn}
\pp_{\ell,\sigma,\tau}((\sqrt{s_i s_i ^\prime })_{1\leq i\leq \ell}  )^2<
\pp_{\ell,\sigma,\tau}(s)\cdot \pp_{\ell,\sigma,\tau}(s^\prime ).
\end{equation}  
\end{lemma}

\begin{proof}   For simplicity and as in Lemma  \ref{lem: Zspan},  
 we will only prove the lemma for the case when $\tau=\mathrm{Id}$ is the identity permutation.
We write $f := \pp_{\ell,\sigma, \mathrm{Id}}$ on $\R_+^{\ell}$, 
 so that, with the coordinate functions on $\R^{\Sigma(1)}$ arranged as in \eqref{eq: sigmaSigma1}, we have
\begin{equation}\label{eq: fpp}f(s_1,\dots,s_\ell)=
\pp_{\ell,\sigma,\mathrm{Id}}(s_1,\dots,s_\ell)=\pp(s_1,\cdots, s_\ell,0,\cdots,0,1\cdots,1)
\end{equation}
for $s_1,\dots,s_\ell\in \R_+^{\ell}$.  
As in \eqref{eq: tildef},
	we consider the associated function $f^\sharp:\R^{\ell }\to \R$
		given by $f^\sharp(t_1, \ldots, t_{\ell }) := f(e^{t_1}, \ldots, e^{t_{\ell}})$.
By Lemma \ref{lem: Zspan}(i), the Newton polytope 
	$\Phi(f)$ of $f$
		has affine dimension $\ell $.
Upon regarding $f$ as a polynomial on $\C^\ell$, it also follows from \eqref{eq:  pinequalitycn}
	that $|f(z_1, \ldots, z_{\ell })| \le f(|z_1|, \ldots, |z_{\ell }|)$
		for all $(z_1, \ldots, z_{\ell }) \in \C^{\ell }$.			
Hence, by Lemma \ref{lem: positiveHessianInInterior} and \eqref{eq: HessJlogf},
	the Hessian matrix 
	\begin{equation}\label{eq: hessianpositive}
	\Big( \frac{\partial^2}{\partial t_i \partial t_j} \log {f}^\sharp(t)\Big)_{1 \le i, j \le \ell }			
		\text{ is positive definite for all }t \in \R^{\ell}.
		\end{equation}
			   Let $(s_1^*,\cdots, s_k^*)\in (\R_+^k)^\circ$, and let
				$t_i^*=\log s_i^*$ for each $1\leq i\leq k$.  Let $g:[0,s_k^*]\to\R$ be the function given by $ g(s)=\log f(s_1^*,\cdots, s_{k-1}^*,s)$ for $s\in [0,s_k^*]$, and let $h:(-\infty,t_k^*]\to\R$
				be given by $h(t)=\log f^\sharp(t_1^*,\cdots, t_{k-1}^*,t)$ for $t\in (-\infty,t_k^*]$.
Then it follows readily from \eqref{eq: hessianpositive} that $h^{\prime\prime}(t)>0$ for all $t\in (-\infty,t_k^*]$, and thus $h^{\prime}(t)$ is a strictly increasing function in $t$ for $t\in (-\infty,t_k^*]$.  
By the chain rule, one easily check that $h^\prime(t)=g^\prime(e^t)\cdot e^t$, and thus it follows that
the function $\mu(s):=s\cdot g^\prime(s)$ is a strictly increasing function in $s$ for $s\in (0,s_k^*]$.  Together with continuity of $\mu$ on the interval $[0, s_k^*]$ and the fact that $\mu(0)=0$, it follows that 
$\mu(s)>0$ for all $s\in (0, s_k^*]$, and thus $g^\prime(s)>0$ for all $s\in (0, s_k^*]$.  Hence $g(s)$
(and thus also $e^{g(s)}=f(s_1^*,\cdots, s_{k-1}^*,s)$) is a strictly increasing function in $s$ for 
$s\in (0, s_k^*]$.  Together with \eqref{eq: fpp} (and upon renaming each $s_i^*$ as $s_i$), one obtains \eqref{eq: plincreasing} readily, and this finishes the proof of (i).  We proceed to prove (ii).  First from \eqref{eq: hessianpositive}, 
one knows that $\log {{f}^\sharp}$ is a strictly convex function on $\R^{n}$
				(see e.g. \cite[p. 37]{BL06}).   
In particular,
	we have 
\begin{equation}\label{eq: convexfsharp}  
\log {f}^\sharp( \frac{t + t^\prime}{2}) < \frac{1}{2}(\log{f}^\sharp(t) +\log{f}^\sharp(t^\prime))
\end{equation}  
for all $t=(t_1,\dots,t_{n})\neq  t^\prime=(t_1^\prime,\dots,t_{n}^\prime) \in \R^{n }$.
Upon exponentiating both sides of \eqref{eq: convexfsharp} and letting $t_i=\log s_{i},~t_i^\prime=\log s_{i}^\prime$ for each $i$, one obtains \eqref{eqlem: convexEqn} readily.  
\end{proof}

\begin{proposition} \label{prop: pPSGCS}
Suppose $\pp$ satisfies (Pos1), (Pos2) and (Pos3).  Then the following statements hold:
\par\noindent
(i) For all $x=(x_\rho)_{\rho\in\Sigma(1)}, ~y=(y_\rho)_{\rho\in\Sigma(1)}
\in \R_+^{\Sigma(1)}\setminus
(Z(\Sigma)\cap\R_+^{\Sigma(1)})$ such that $(\sqrt{x_\rho y_\rho})_{\rho\in\Sigma(1)} 
\in  \R_+^{\Sigma(1)}\setminus
(Z(\Sigma)\cap\R_+^{\Sigma(1)})$, we have
\begin{equation}\label{eq: pxy2pxpy}
\pp((\sqrt{x_\rho y_\rho})_{\rho\in\Sigma(1)} )^2 \le \pp(x)\cdot \pp(y).
\end{equation}
Furthermore, the equality in \eqref{eq: pxy2pxpy} holds if and only if there exists $g\in G\cap (\R_+^{\Sigma(1)})^\circ $ such that
$g\cdot x=y$.
\par\noindent
(ii)  For all $z,w\in\C^{\Sigma(1)}\setminus Z(\Sigma)$ such that $\Xi(z)\neq \Xi(w)$, one has
\begin{equation}\label{eq: PtildeSCGS}
|\PP(z,\overline{w})|^2 < \PP(z,\overline{z})\cdot \PP(w,\overline{w}).
\end{equation}
\par\noindent
(iii)
$\widehat{P}$ satisfies the SGCS inequality.
\end{proposition}

\begin{proof}  Let $x=(x_\rho)_{\rho\in\Sigma(1)}, ~y=(y_\rho)_{\rho\in\Sigma(1)}
\in \R_+^{\Sigma(1)}\setminus
(Z(\Sigma)\cap\R_+^{\Sigma(1)})$ be such that $(\sqrt{x_\rho y_\rho})_{\rho\in\Sigma(1)} 
\in  \R_+^{\Sigma(1)}\setminus
(Z(\Sigma)\cap\R_+^{\Sigma(1)})$.  To prove (i), we first establish the equality 
 in \eqref{eq: pxy2pxpy} under the assumption that 
 there exists $g=(t_{\rho})_{\rho\in\Sigma(1)}\in G\cap (\R_+^{\Sigma(1)})^\circ $ such that
$g\cdot x=y$.  Let $g_\sharp:=(\sqrt{t_{\rho}})_{\rho\in\Sigma(1)}$.  Then one easily sees that 
$g_\sharp\in G\cap (\R_+^{\Sigma(1)})^\circ $,  and
\begin{equation}\label{eq: g2ppxy}
g_\sharp\cdot x=g_\sharp^{-1}\cdot y=  (\sqrt{x_{\rho}y_{\rho}})_{\rho\in\Sigma(1)}).
\end{equation}
Let $\chi^L$ be the group character in \eqref{eq: Gcharacter}.   
Together with \eqref{eq: pfunctionaleqn}, one has
\begin{align}\notag
\pp(\sqrt{x_{\rho}y_{\rho}})_{\rho\in\Sigma(1)})^2&=\pp (g_\sharp\cdot x)\cdot\pp(g_\sharp^{-1}\cdot y)
=\chi^L (g_\sharp)\cdot \pp(x)\cdot \chi^L (g_\sharp^{-1})\cdot \pp(y)\\
\notag
& =\chi^L (g_\sharp)\cdot \pp(x)\cdot \chi^L (g_\sharp)^{-1}\cdot \pp(y)  =\pp(x)\cdot \pp(y).
\end{align}
To complete the proof (i), it remains to establish the strict inequality in \eqref{eq: pxy2pxpy} for the case when $y\notin (G\cap (\R_+^{\Sigma(1)})^\circ )\cdot x$.  Since $(\sqrt{x_\rho y_\rho})_{\rho\in\Sigma(1)} 
\in  \R_+^{\Sigma(1)}\setminus
(Z(\Sigma)\cap\R_+^{\Sigma(1)})$, there exists $\sigma\in\Sigma(n)$ such that 
$\sqrt{x_\rho y_\rho}>0$ (and thus $x_\rho>0$ and $y_\rho>0$) for all $\rho\in \Sigma(1)\setminus\sigma(1)$.  
We arrange the coordinate functions on $\R^{\Sigma(1)}$ as in \eqref{eq: sigmaSigma1}.
Then it follows readily from Lemma \ref{lem: gx} that
there exist $g=(t_{i})_{1\leq i\leq |\Sigma(1)|},\, g^\prime=(t_{i}^\prime)_{1\leq i\leq |\Sigma(1)|}\in G\cap (\R_+^{\Sigma(1)})^\circ   $ and $s=(s_{i})_{1\leq i\leq n},  s^\prime=(s^\prime_{i})_{1\leq i\leq n}\in  (\R_+^{\sigma(1)})^\circ$
such that
\begin{equation} \label{eq: gxgprime}
g\cdot x=(s,1,\cdots,1)\quad\text{and}\quad g^\prime\cdot y=(s^\prime,1,\cdots,1).
\end{equation}  
Since
$y\notin (G\cap (\R_+^{\Sigma(1)})^\circ )\cdot x$, it follows readily that $s\neq s^\prime$.  
Let $g^{\prime\prime}:=(\sqrt{t_{i} t_{i}^\prime})_{1\leq i\leq |\Sigma(1)|}$.  Then
one easily sees from \eqref{eq:  defG} and \eqref{eq: gxgprime} that $g^{\prime\prime}\in G\cap (\R_+^{\Sigma(1)})^\circ$ and 
\begin{align}\label{eq: gppsqrt}
g^{\prime\prime}\cdot (\sqrt{x_{i} y_{i}^\prime})_{1\leq i\leq |\Sigma(1)|})
&=(\sqrt{s\cdot s^\prime} , 1,\cdots,1),\quad\text{where}\\
\sqrt{s\cdot s^\prime}:&=(\sqrt{s_{1} s_{1}^\prime}, \ldots, \sqrt{s_{n} s_{n}^\prime}).
\end{align}
Let $\ell$ be the number of non-zero entries of $\sqrt{s\cdot s^\prime}$ (so that $ 0\leq\ell \leq n$.
First we consider the case when $\ell\geq 1$.  We further rearrange the first $n$ coordinate functions of 
$\R^{\Sigma(1)}$, so that we may assume that
\begin{equation}\label{eq: sisiprime}
\sqrt{s_{i} s_{i}^\prime}>0 \quad\text{for }1\leq i\leq\ell\quad\text{and}\quad 
\sqrt{s_{i} s_{i}^\prime}=0\quad\text{for }\ell+1\leq i\leq n
\end{equation}
(the last equality is ignored if $\ell=n$).  Now we let
\begin{equation}\label{eq: hatsisiprime}
\widehat{s}:=(s_i)_{1\leq i\leq\ell}, \quad \widehat{s^\prime}:=(s_i^\prime)_{1\leq i\leq\ell}
\quad\text{and}\quad \sqrt{\widehat{s\cdot s^\prime}}:=(\sqrt{s_i\cdot s_i^\prime})_{1\leq i\leq\ell}.
\end{equation}
Then it is easy to see that 
 $\displaystyle \widehat{s}, \,\widehat{s^\prime}, \,\sqrt{\widehat{s\cdot s^\prime}}\in  (\R_+^{\ell})^\circ$.
If $\widehat{s}\neq\widehat{s^\prime}$, then it follows from Lemma \ref{lem: pconvex}(ii) that
\begin{equation}\label{eq: firststrictinequality}
\displaystyle\pp_{\ell,\sigma,\text{Id}}(\sqrt{\widehat{s\cdot s^\prime}})^2<\pp_{\ell,\sigma,\text{Id}}(\widehat{s})\cdot 
\pp_{\ell,\sigma,\text{Id}}(\widehat{s^\prime})\leq \pp_{n,\sigma,\text{Id}}(s)\cdot \pp_{n,\sigma,\text{Id}}(s^\prime),
\end{equation}
where the last inequality follows from using 
Lemma \ref{lem: pconvex}(i) repeatedly if necessary.
On the other hand, if $\displaystyle\widehat{s}=\widehat{s^\prime}$ (so that $\displaystyle\sqrt{\widehat{s\cdot s^\prime}}=\widehat{s}=\widehat{s^\prime}$), then since $s\neq s^\prime$, it follows readily that either $s_i>0$ for some $\ell+1\leq i\leq n$, or $s_i^\prime>0$ for some $\ell+1\leq i\leq n$.  By using 
Lemma \ref{lem: pconvex}(i) repeatedly, one sees that either 
$\pp_{\ell,\sigma,\text{Id}}(\widehat{s})<\pp_{n,\sigma,\text{Id}}(s)$ or 
$\pp_{\ell,\sigma,\text{Id}}(\widehat{s^\prime})<\pp_{n,\sigma,\text{Id}}(s^\prime)$ respectively.
Thus in this case, we have
\begin{equation}\label{eq: secondtrictinequality}
\displaystyle\pp_{\ell,\sigma,\text{Id}}(\sqrt{\widehat{s\cdot s^\prime}})^2=\pp_{\ell,\sigma,\text{Id}}(\widehat{s})\cdot 
\pp_{\ell,\sigma,\text{Id}}(\widehat{s^\prime})<\pp_{n,\sigma,\text{Id}}(s)\cdot \pp_{n,\sigma,\text{Id}}(s^\prime) .
\end{equation}
Note that from \eqref{eq: sisiprime}, one has $\displaystyle\pp_{\ell,\sigma,\text{Id}}(\sqrt{\widehat{s\cdot s^\prime}})=\pp_{n,\sigma,\text{Id}}(\sqrt{s\cdot s^\prime})$.  
Together with \eqref{eq: firststrictinequality}, it follows that when $\ell>0$,
we have 
\begin{align}\label{eq: 3rdstrictinequality}
\displaystyle \pp_{n,\sigma,\text{Id}}(\sqrt{s\cdot s^\prime})^2&<\pp_{n,\sigma,\text{Id}}(s)\cdot \pp_{n,\sigma,\text{Id}}(s^\prime) ,\quad \text{or equivalently},\\
\label{eq: 4rdstrictinequality}
\pp(\sqrt{s\cdot s^\prime} , 1,\cdots,1)^2&< \pp(s, 1,\cdots,1)\cdot\pp(s^\prime, 1,\cdots,1))
\end{align}
for both cases when $\widehat{s}\neq\widehat{s^\prime}$ and when $\widehat{s}=\widehat{s^\prime}$.
We remark that in the case when $\ell=0$, one sees that \eqref{eq: 4rdstrictinequality} still holds by using an argument similar to the inequality in \eqref{eq: secondtrictinequality}.
From \eqref{eq: pfunctionaleqn}, \eqref{eq: gxgprime} and \eqref{eq: gppsqrt}, one sees that
\begin{align}\label{eq: chixy1}
\pp(s,1,\cdots,1)&=\pp(g\cdot x)=\chi^L(g)\cdot
\pp(x) ,
\\
\label{eq: chixy2}
\pp(s^\prime,1,\cdots,1)&=\pp(g^\prime\cdot y)=\chi^L(g^\prime)\cdot
\pp(y) ,\quad\text{and}
\\
\label{eq: chixy3}
\pp(\sqrt{s\cdot s^\prime} , 1,\cdots,1)
&=\chi^L(g^{\prime\prime})\cdot \pp((\sqrt{x_\rho y_\rho})_{\rho\in\Sigma(1)} )
.
\end{align}
Since $(g^{\prime\prime})^2=g\cdot g^\prime\in G\cap (\R_+^{\Sigma(1)})^\circ$,
 it follows readily that
\begin{equation}\label{eq: chiggprime}
(\chi^L(g^{\prime\prime}))^2=\chi^L(g)\cdot \chi^L(g^\prime)>0.
\end{equation}
By combining \eqref{eq: 4rdstrictinequality}, \eqref{eq: chiggprime},
\eqref{eq: chixy1}, \eqref{eq: chixy2}, \eqref{eq: chixy3} and \eqref{eq: chiggprime}, one sees that
\begin{equation}\label{eq: ppxyinequality}
\pp((\sqrt{x_\rho y_\rho})_{\rho\in\Sigma(1)} )^2
< \pp(x)\cdot \pp(y),
 \end{equation}
and we have finished proof of (i).
 We proceed to prove (ii).
 Let  $z=(z_\rho)_{\rho\in\Sigma(1)},\,w=  (w_\rho)_{\rho\in\Sigma(1)}\in\C^{\Sigma(1)}\setminus Z(\Sigma)$ be such that $\Xi(z)\neq \Xi(w)$.  For notational convenience, we denote
 $z\cdot \overline{w}:=(z_\rho\cdot\overline{w_\rho})_{\rho\in\Sigma(1)}$, $z\cdot\overline{z}=(|z_\rho|^2)_{\rho\in\Sigma(1)}$ and $w\cdot\overline{w}:=(|w_\rho|^2)_{\rho\in\Sigma(1)}$.
 Then from \eqref{eq: Phihomogenizationofp} and \eqref{eq: PPhihomogenizationofp}, one sees that
$
\PP(z,\overline{w})=\pp(z\cdot \overline{w})$.  
We consider the following two cases:
\par\noindent
$\underline{Case~(a)}:~$  When $z\cdot \overline{w}\notin G\cap U(1)^{\Sigma(1)})\cdot \R_+^{\Sigma(1)}$, it follows from (Pos3) that 
\begin{equation}\label{eq:  ppzw1}
|\pp(z\cdot \overline{w})|^2<|\pp((|z_\rho||w_\rho|)_{\rho\in\Sigma(1)})|^2
\leq \pp(z\cdot \overline{z})\cdot\pp(w\cdot \overline{w}),
\end{equation}
where the last inequality follows from (i).
\par\noindent
$\underline{Case~(b)}$:~When $z\cdot \overline{w}\in G\cap U(1)^{\Sigma(1)})\cdot \R_+^{\Sigma(1)}$, one easily sees that we have  $z\cdot\overline{w}=(e^{i\theta_\rho})_{\rho\in\Sigma(1)}\cdot
 (|z_\rho||w_\rho|)_{\rho\in\Sigma(1)})$, where $(e^{i\theta_\rho})_{\rho\in\Sigma(1)}\in G\cap U(1)^{\Sigma(1)})$.  First we see from \eqref{eq: Gcharacter} that
 \begin{equation}\label{eq: functionaleitheta}
| \pp(z\cdot\overline{w})|=|\chi^L((e^{i\theta_\rho})_{\rho\in\Sigma(1)})|\cdot |\pp ( (|z_\rho||w_\rho|)_{\rho\in\Sigma(1)}))|=|\pp ( (|z_\rho||w_\rho|)_{\rho\in\Sigma(1)}))|,
 \end{equation}
 where the last equality follows from the equality $|\chi^L((e^{i\theta_\rho})_{\rho\in\Sigma(1)})|=1$, which can be seen easily from \eqref{eq: Gcharacter}.
Suppose that $ (|z_\rho|)_{\rho\in\Sigma(1)}=g\cdot
 (|w_\rho|)_{\rho\in\Sigma(1)}$ for some $g\in G\cap(\R_+^{\Sigma(1)})^\circ$.  
By considering the components individually, one easily sees that $z=g\cdot (e^{i\theta_\rho})_{\rho\in\Sigma(1)}
\cdot w$ with $g\cdot (e^{i\theta_\rho})_{\rho\in\Sigma(1)}\in G$, contradicting the assumption that 
$\Xi(z)\neq\Xi(w)$.   Thus we must have $(|z_\rho|)_{\rho\in\Sigma(1)}\notin 
(G\cap(\R_+^{\Sigma(1)})^\circ)\cdot (|w_\rho|)_{\rho\in\Sigma(1)}$, and then it follows from (i) that 
\begin{equation}\label{eq:  ppzw2}|\pp(z\cdot \overline{w})|^2=
|\pp((|z_\rho||w_\rho|)_{\rho\in\Sigma(1)})|^2
<\pp(z\cdot \overline{z})\cdot\pp(w\cdot \overline{w}),
\end{equation}
where the first equality follows from \eqref{eq: functionaleitheta}.  
\par\noindent
Thus in both cases, we have $|\pp(z\cdot \overline{w})|^2< \pp(z\cdot \overline{z})\cdot\pp(w\cdot \overline{w})$, which upon rewritten in terms of $\PP$,  leads to \eqref{eq: PtildeSCGS}, and we have finished the proof of (ii).  Finally from the lifting $\Xi^*L\cong (\C^{\Sigma(1)}\setminus Z(\Sigma))\times 
\C$ in Remark \ref{rem: triviallb}, one sees
from \eqref{eq: Pdef2} and \eqref{eq: PPhihomogenizationofp}
that $\widehat{P}$ lifts to the function $\PP$ on 
$(\C^{\Sigma(1)}\setminus Z(\Sigma))\times \overline{(\C^{\Sigma(1)}\setminus Z(\Sigma))}$.  
Then  the inequality in \eqref{eq: PtildeSCGS} leads readily to the SGCS inequality for $\widehat{P}$ (cf. \eqref{eq: SGCSineq}), and this finishes the proof of (iii).
\end{proof}

We conclude this section with the following

\begin{proof}[Proof of Proposition \ref{prop: P3properties}]
Proposition \ref{prop: P3properties} follows directly from
Proposition \ref{prop: Pos1Pos3impliesPpd},
Proposition \ref{prop: StronglyPseudoconvex} and
Proposition \ref{prop: pPSGCS}(iii).
\end{proof}

\section{Proof of Theorem \ref{thm: mainThm}} \label{sec: proofOfTheorem}

\begin{lemma} \label{lem: posCoeffs}
Let $p \in \R[x_1^{\pm 1},\dots,x_n^{\pm 1}]$ be a Laurent polynomial whose Newton polytope $\Phi$ is an
$n$-dimensional smooth lattice polytope. Let $\pp$ be the $\Phi$-homogenization of $p$, and let $\Sigma$ be the associated normal fan of $\Phi$.  Suppose $p$ has fully positive coefficients.
Then 
\begin{equation} \label{eq: positiveOnOrthant}
\pp(x) > 0\quad \text{for all } x \in \R_+^{\Sigma(1)} \setminus (Z(\Sigma)\cap \R_+^{\Sigma(1)} )
,
\end{equation}	
and $\pp$ satisfies (Pos1), (Pos2) and (Pos3).
\end{lemma}

\begin{proof} Write $\displaystyle p=\sum_{m\in\Phi\cap M}c_mx^m$.  Since $\pp$ has fully positive coefficients, we have
$c_m>0$ for all $m\in\Phi\cap M$.  
Write $\displaystyle
\pp:=\sum_{m\in \Phi\cap M }c_m\cdot \prod_{\rho\in\Sigma(1)}x_\rho^{\langle m,u_\rho\rangle+
a_\rho}$ as in \eqref{eq: Phihomogenizationofp}.
Let $x=(x_\rho)_{\rho\in \Sigma(1)} \in \R_+^{\Sigma(1)} \setminus (Z(\Sigma)\cap \R_+^{\Sigma(1)} )$. 
Then from \eqref{eq:  irrelevantset}, there exists $\sigma\in\Sigma(n)$ such that 
$x_\rho>0$ for all $\rho\in\Sigma(1)\setminus\sigma(1)$.  Let $v\in\Phi$ be the vertex of $\Phi$ that corresponds to $\sigma$ under the Orbit-Cone Correspondence.  Then 
for all $\rho\in\sigma(1)$, we have $\langle v,u_\rho\rangle+a_\rho=0$.  It follows that 
\begin{equation}\label{eq: pppositive}
\pp(x)\geq c_v\cdot \prod_{\rho\in \Sigma(1)\setminus\sigma(1)} x_\rho^{\langle v,u_\rho\rangle+a_\rho}>0,
\end{equation}
which gives \eqref{eq: positiveOnOrthant}.  For each $\sigma\in \Sigma(n)$, one easily sees from 
\eqref{eq: defesigma} and \eqref{eq:  irrelevantset} that $e^{(\sigma)}\in \R_+^{\Sigma(1)} \setminus (Z(\Sigma)\cap \R_+^{\Sigma(1)} )$, and thus it follows from \eqref{eq: positiveOnOrthant} that
$\pp(e^{(\sigma)})>0$.  Thus $\pp$ satisfies (Pos1).  To show that $\pp$ satisfies (Pos2),  we let $x=(x_\rho)_{\rho\in\Sigma(1)}
\in F_{\rho_\circ}(\R_+^{\Sigma(1)}) \setminus
  											 (Z(\Sigma)\cap F_{\rho_\circ}(\R_+^{\Sigma(1)}) )$ with $\rho_o\in
  											 \Sigma(1)$, so that $x_{\rho_o}=0$.
  											 As above, since $x\notin Z(\Sigma)$, it follows that there
  											 exists $\sigma\in\Sigma(n)$ such that $x_\rho>0$ for all
  											 $\rho\in\Sigma(1)\setminus\sigma(1)$.  Again as above, let $v\in\Phi$ be the vertex of $\Phi$ that corresponds to $\sigma$ under the Orbit-Cone Correspondence. Since $\Phi$ is an $n$-dimensional smooth lattice polytope, there exist $n$ edges $\{E_\rho\}_{\rho\in \sigma(1)}$ of $\Phi$ (with $\displaystyle E_\rho=\Phi\cap\bigcap_{\rho^\prime\in\sigma(1)\setminus\{\rho\}}
H_{\rho^\prime}$) emanating from $v$, each containing a unique lattice point $w_{\rho}(\in\Phi\cap M)$ adjacent to $v$, and such that $\{w_\rho-v\}_{\rho\in\sigma(1)}$ forms a $\Z$-basis of $M$.  
Recall that $\sigma$ is an $n$-dimensional smooth cone in $N_\R$.
It is easy to check that $\{u_\rho\}_{\rho\in \sigma(1)}$ forms a $Z$-basis of $N$ dual to $\{w_\rho-v\}_{\rho\in\sigma(1)}$, and it follows readily that
$\langle  w_{\rho},u_{\rho^\prime}\rangle+a_{\rho^\prime}=\delta_{\rho\rho^\prime}$ for $\rho,\rho^\prime\in \sigma(1)$.  Here $\delta_{\rho\rho^\prime}$ denotes the Kronecker delta.  Then it follows readily that
  			\begin{equation}\label{eq: partialpppositive}
\frac{\partial \pp}{ \partial x_{\rho_\circ}} (x) \geq c_{w_\rho }\cdot \prod_{\rho\in \Sigma(1)\setminus\sigma(1)} x_\rho^{\langle w_{\rho_\circ},u_\rho\rangle+a_\rho}
> 0,
\end{equation}								 
which implies that $\pp$ satisfies (Pos2).  We proceed to show that $\pp$ satisfies (Pos3).  Let $z=(z_\rho)_{\rho\in\Sigma(1)} \in \C^{\Sigma(1)} \setminus (Z(\Sigma)\cup (G\cap U(1)^{\Sigma(1)})\cdot \R_+^{\Sigma(1)})$.  
Since $z\notin Z(\Sigma)$, there exists $\sigma\in\Sigma(n)$ such that 
$z_\rho\neq 0$ {for all }$\rho\in\Sigma(1)\setminus\sigma(1)$.  Then as in the proof of Proposition \ref{prop: pPSGCS}, by using an element of $G$ as given in Remark \ref{rem: gx111} and using the functional
equation in \eqref{eq: pfunctionaleqn}, it suffices to consider the case when 
\begin{equation}\label{eq: 111}
z_\rho=1\text{  for all }\rho\in\Sigma(1)\setminus \sigma(1),\text{ so that }B:=\{\rho\in\Sigma(1)\,\big|\,z_\rho=0\}\subset \sigma(1).
\end{equation}
For each $\rho\in \sigma(1)\setminus B$, we fix $\theta_\rho\in\R$ such that $z_\rho=e^{i\theta_\rho}\cdot |z_\rho|$.  For each $\rho\in B\cup (\Sigma(1)\setminus \sigma(1))$, we simply let $\theta_\rho=0$.
Then one easily sees that
 \begin{equation}\label{eq: zeitheta}
z=(e^{i\theta_\rho})_{\rho\in\Sigma(1)}\cdot
(|z_\rho|)_{\rho\in \Sigma(1)}\quad\text{with}\quad (e^{i\theta_\rho})_{\rho\in\Sigma(1)}\in U(1)^{\Sigma(1)}.
\end{equation}
Let
\begin{equation}\label{eq: Psisubset}
\Psi:=\Phi\cap M\cap\bigcap_{\rho\in B}H_\rho,
\end{equation}
where $H_\rho$ is as in \eqref{eq:defHj}.  Then upon ignoring the zero terms of $\pp(z)$ arising
from those $z_\rho$'s with $\rho\in B$ and ignoring the trivial factors of each monomial term of $\pp(z)$
(arising from the exponent being zero or the condition that $z_\rho=1$ for $\rho\in\Sigma(1)\setminus \sigma(1))$,
one sees that
\begin{align}\label{eq: ppdifferent}
\pp(z)&=
\sum_{m\in\Psi}c_m\cdot\prod_{\rho\in\sigma(1)\setminus B}z_\rho^{\langle m,u_\rho\rangle+a_\rho},
\quad\text{so that}\\
\label{eq: ppleq}
|\pp(z)| &=\big|\sum_{m\in \Psi}  c_m\cdot\prod_{\rho\in\sigma(1)\setminus B}z_\rho^{\langle m,u_\rho\rangle+a_\rho}
\big|\\
\notag
&
\leq 
\sum_{m\in \Psi}  c_m\prod_{\rho\in\sigma(1)\setminus B}|z_\rho|^{\langle m,u_\rho\rangle+a_\rho}
=\pp((|z_\rho|)_{\rho\in \Sigma(1)}).
\end{align}  Thus it remains to show that if the non-trivial equality in \eqref{eq: ppleq} holds, then 
$z\in ( G\cap U(1)^{\Sigma(1)})\cdot \R_+^{\Sigma(1)}$, or equivalently, $(e^{i\theta_\rho})_{\rho\in\Sigma(1)}\in G$.  To see this, we suppose that the non-trivial equality  in \eqref{eq: ppleq} holds.  This implies readily that the arguments of the (non-zero) complex numbers $\displaystyle \prod_{\rho\in\sigma(1)\setminus B}z_\rho^{\langle m,u_\rho\rangle+a_\rho}$ are the same for all $m\in\Psi$, so that there exists $\kappa\in\R$ such that 
\begin{equation}\label{eq: sumtheta}
\sum_{\rho\in\sigma(1)\setminus B}(\langle m,u_\rho\rangle+a_\rho)\theta_\rho
\equiv\kappa\text{   mod }2\pi  \quad
\text{for all }m\in\Psi.
\end{equation}
Let $v\in\Phi$ be the vertex that corresponds to $\sigma$ under the Orbit-Cone Correspondence, so that $v=\bigcap_{\rho\in\sigma(1)} H_\rho$.  Then from \eqref{eq: Psisubset}, one easily sees that $v\in\Psi$.  
Let $\{w_\lambda\}_{\lambda\in\sigma(1)}$ be the unique lattice points of $\Phi$ adjacent to $v$ along the edges $\{E_\lambda\}_{\lambda\in\sigma(1)}$ 
emanating from $v$ as in \eqref{eq: partialpppositive}.  Note that 
\begin{equation}\label{eq: wlambda}
w_\lambda\in \Phi\cap M\cap \bigcap_{\rho\in \sigma(1)\setminus \{\lambda\} }H_{\rho}\quad\text{for each }\lambda\in\sigma(1).
\end{equation}
 Then similar to $v$, it is easy to see from \eqref{eq: Psisubset} and \eqref{eq: wlambda} that $w_\lambda\in\Psi$ for each $\lambda\in \sigma(1)\setminus B$.  Thus the equation in \eqref{eq: sumtheta} is satisfied by $m=v$ and $m=w_\lambda$ for $\lambda\in \sigma(1)\setminus B$.  
 Upon
subtracting these two equations and using the fact that $\theta_\rho=0$ for $\rho\in B$, one gets
 \begin{equation}\label{eq: sumthetaminusv}
\sum_{\rho\in\sigma(1)}\langle w_\lambda-v,u_\rho\rangle\theta_\rho
\equiv 0\text{   mod }2\pi  \quad\text{for all }\lambda \in \sigma(1)\setminus B.
\end{equation}
 Now we consider the case when $\lambda\in B$ (so that $\sigma(1)\setminus B\subset 
 \sigma(1)\setminus\{\lambda\}$).  In this case, it follows from \eqref{eq: wlambda} that for all $\rho\in \sigma(1)\setminus B$, one has
$\langle w_\lambda,u_\rho\rangle+a_\rho=0$ and thus $\langle w_\lambda-v,u_\rho\rangle=0$.  Together with the fact that $\theta_\rho=0$ for all $\rho\in B$. it follows that we have
 \begin{equation}\label{eq: sumthetaminusv2}
\sum_{\rho\in\sigma(1)}\langle w_\lambda-v,u_\rho\rangle\theta_\rho
=0 \quad\text{for all }\lambda \in B.
\end{equation}
Since $\Phi$ is an $n$-dimensional smooth lattice polytope, it follows that $\{w_\rho-v\}_{\rho\in\sigma(1)}$ forms a $\Z$-basis of $M$.  Together with \eqref{eq: sumthetaminusv} and \eqref{eq: sumthetaminusv2}, it follows readily that 
 \begin{equation}\label{eq: sumthetaminusvm}
\langle m,\sum_{\rho\in\sigma(1)}\theta_\rho u_\rho\rangle
\equiv 0   \text{   mod }2\pi \quad\text{for all }m\in M.
\end{equation}
Upon exponentiating \eqref{eq: sumthetaminusvm} and using the fact that $\theta_\rho$ for all $\rho\in \Sigma(1)\setminus \sigma(1)$, it follows that
 \begin{equation}\label{eq: exponthetaG}
\prod_{\rho\in\Sigma(1)}(e^{i\theta_\rho})^{
\langle m,u_\rho\rangle}
=1\quad\text{for all }m\in M, 
\end{equation}
which together with \eqref{eq:  defG}, implies that $(e^{i\theta_\rho})_{\rho\in\Sigma(1)}\in G$, and 
we have finished the proof that $\pp$ satisfies (Pos3).
\end{proof}

We are ready to give the proof of Theorem \ref{thm: mainThm} as follows:

\begin{proof}[Proof of Theorem \ref{thm: mainThm}]  Let $p \in \R[x_1^{\pm 1},\dots,x_n^{\pm 1}]$ be a Laurent polynomial whose Newton polytope $\Phi$ is an
$n$-dimensional smooth lattice polytope. Let $\pp$ be the $\Phi$-homogenization of $p$.  
\par\noindent
$\underline{\eqref{thm: Condition}\implies\eqref{thm: EventualPos}}$:  
Suppose $\pp$ satisfies (Pos1), (Pos2) and (Pos3).   Let $\widehat{P}$ be the Hermitian algebraic function 
on the line bundle $L$ as in \eqref{eq: Pdef2}, where $L$ as in \eqref{eq:  defL}.  
 Then from Proposition \ref{prop: P3properties} and 
Theorem \ref{TheoremCD99} (with $R=\widehat{P}$, $E=\mathcal O_{X_\Sigma}$ and $Q=1$), one knows that  there exists $k_o > 0$ such that for each integer $k \ge k_o$,
            the tensor power $\widehat{P}^k $ (as a Hermitian algebraic function on $L^{\otimes k}$) is a maximal sum of Hermitian squares.
            It is easy to see that $\widehat{P}^k $ is the Hermitian algebraic function
            on $L^{\otimes k}$ induced from the Laurent polynomial $p^k$.  Note also that
            the Newton polytope $\Phi(p^k)$ of $p^k$ is simply $k\cdot \Phi:=\{k\cdot x\,\big|\, x\in \Phi\}$, which is also an $n$-dimensional smooth lattice polytope.  
            Then from Proposition  \ref{prop: AssocBihomIsSOS},
            it follows that $p^k$ has fully coefficients for each such $k$.
            \par\noindent
$\underline{\eqref{thm: EventualPos} \implies \eqref{thm: Condition}}$: 
Suppose that there exists $k_o > 0$ such that for each integer $k \ge  k_o$, 
								\, $p^k$ has fully positive coefficients.  
								From Lemma \ref{lem: posCoeffs} (with $\pp$ there replaced by $\pp^k$),
	one knows that $\pp(x)^k>0$ for all $ x \in \R_+^{\Sigma(1)} \setminus (Z(\Sigma)\cap \R_+^{\Sigma(1)} ) $, and $\pp^k$ satisfies (Pos1), (Pos2) and (Pos3).  
By choosing an odd $k\geq k_o$ and then taking the $k$-th root, one immediately sees that $\pp(x)>0$ for all $ x \in \R_+^{\Sigma(1)} \setminus (Z(\Sigma)\cap \R_+^{\Sigma(1)} )$, and 
$\pp$ also satisfies (Pos1) and (Pos3).
Since $\pp(x)>0$ for all $ x \in \R_+^{\Sigma(1)}\setminus (Z(\Sigma)\cap \R_+^{\Sigma(1)} )$ and we have
\begin{equation}
\frac{\partial (\pp^k)}{\partial x_\rho}(x) = k \cdot \pp(x)^{k - 1} \cdot \frac{\partial \pp}{\partial x_\rho}(x) \quad\text{for each }
\rho\in \Sigma(1),
\end{equation}
it follows readily from  (Pos2) for  $\pp^k$ that
$\pp$ also satisfies (Pos2).		
\end{proof}

\section{Application to polynomial spectral radius functions} \label{sec: Application}

In this section, we apply Theorem \ref{thm: mainThm} to 
deduce Corollary \ref{cor: betaApplication}.  As the proof is similar to \cite[Section 5]{TT18}, 
 we refer the reader to \cite[Section 5]{TT18} for a convenient recollection of the background results needed as well as the definitions of the various terms in the corollary.

\begin{proof}[Proof of Corollary \ref{cor: betaApplication}]
Let $p \in \Z[x_1,\dots,x_n]$ be a polynomial whose Newton polytope $\Phi$ is an $n$-dimensional smooth lattice polytope, and such that the $\Phi$-homogenization $\pp\in \Z[ z_\rho\,\big|\,\rho\in\Sigma(1)]$ of $p$ satisfies (Pos1) and (Pos2).
\par\noindent
$\underline{\eqref{cor: irr}  \implies \eqref{cor: pos3} ~(\text{resp. }
\eqref{cor: aper} \implies \eqref{cor: pos3} } $):
Suppose that there exists an irreducible
(resp. aperiodic) square matrix $A$ over $\Z_+[ z_\rho\,\big|\,\rho\in\Sigma(1)]$ such that $\pp = \beta_A$.
With notation as in Theorem \ref{thm: mainThm}, we let $z=(z_\rho)_{\rho\in\Sigma(1)} \in \C^{\Sigma(1)} \setminus (Z(\Sigma)\cup (G\cap U(1)^{\Sigma(1)})\cdot \R_+^{\Sigma(1)})$.
Since $z\notin Z(\Sigma)$, there exists $\sigma\in\Sigma(n)$ such that 
$z_\rho\neq 0$ {for all }$\rho\in\Sigma(1)\setminus\sigma(1)$.  As in the proof of Proposition \ref{prop: pPSGCS}, by using an element of $G$ as given in Remark \ref{rem: gx111} and using the functional
equation in \eqref{eq: pfunctionaleqn}, it suffices to consider the case when 
$
z_\rho=1\text{  for all }\rho\in\Sigma(1)\setminus \sigma(1)$.   Note that in this case, we have
$(z_\rho)_{\rho\in\sigma(1)}\in \C^{\sigma(1)} \setminus  \R_+^{\sigma(1)}$.
Let $\pp_{n ,\sigma, \mathrm{Id}}$ be as in \eqref{eq: definepmsigma},
	so that
	\begin{equation} \label{eq: formulaForP}
	\pp_{n ,\sigma, \mathrm{Id}}((z_\rho)_{\rho\in\sigma(1)}) = \pp((z_\rho)_{\rho\in\sigma(1)}, 1,\cdots,1)=\pp(z).
	\end{equation}
Then $\pp_{n ,\sigma, \mathrm{Id}} \in \Z [z_\rho\,\big|\, \rho\in\sigma(1)]$
	and $\pp_{n ,\sigma, \mathrm{Id}} = \beta_B$,
		where $B$ is the matrix over $\Z_+ [z_\rho\,\big|\, \rho\in\sigma(1)]]$
			given by $B((z_\rho)_{\rho\in\sigma(1)}) = A((z_\rho)_{\rho\in\sigma(1)},1,\cdots, 1)=A(z)$.
Furthermore, it follows from Lemma \ref{lem: Zspan}(ii) 
	that $S_{\pp_{n ,\sigma, \mathrm{Id}}}$ generates $\Z^{ n}$ as a $\Z$-module.  
	Thus by \cite[Theorem 6.6]{Deangelis942} (cf. also \cite[Lemma 5.2]{TT18}), we have
\begin{equation} \label{eq: hash}
|\pp(z)|=	\pp_{n ,\sigma, \mathrm{Id}}((z_\rho)_{\rho\in\sigma(1)}) | <	\pp_{n ,\sigma, \mathrm{Id}}((|z_\rho|)_{\rho\in\sigma(1)})=\pp((|z_\rho|)_{\rho\in\Sigma(1)}),
\end{equation}	
where the first and last equality follows from \eqref{eq: formulaForP}.				
Hence $\pp$ satisfies (Pos3).
\par\noindent
$\underline{\eqref{cor: pos3} \implies \eqref{cor: irr}~(\text{resp. }\eqref{cor: pos3} 
\implies \eqref{cor: aper})}$: Suppose that $\pp$ also satisfies (Pos3).
Then by Theorem \ref{thm: mainThm},
	there exists $k > 0$ such that $\pp^k$ has fully positive coefficients, so that $\pp^k\in \Z_+[z_\rho\,\big|\, \rho\in \Sigma(1)]$.  By Lemma \ref{lem: posCoeffs}, one also has $\pp^k(x)>0$
	for all $x\in (\R_+^{ \Sigma(1)})^\circ$. 
It follows that
    the $1 \times 1$ matrix $B := (\pp^k)$ is irreducible (resp. aperiodic)  over $\Z_+[z_\rho\,\big|\, \rho\in \Sigma(1)]$, 
    and $\beta_B = \pp^k$.
Hence by {\cite[Theorem 3.3(i) (resp. Theorem 3.5)]{Deangelis942}} (cf. also \cite[Lemma 5.1]{TT18}),  
    there exists 
        an irreducible (resp. aperiodic) 
            square matrix $A$ over $\Z_+[z_\rho\,\big|\, \rho\in \Sigma(1)]$ 
                such that $\pp = \beta_A$.
\end{proof}

\end{document}